\documentclass[12pt,reqno]{amsart}
\usepackage{amsmath,amsfonts,amssymb,amscd,amsthm,amsbsy, color}
\usepackage[dvips]{graphicx}
\usepackage{subfig}
\pdfoutput=1

\usepackage[backend=bibtex, maxbibnames=99]{biblatex}
\newtheorem{theorem}{Theorem}
\newtheorem{meta-thm}[theorem]{Meta-Theorem}

\newtheorem{proposition}[theorem]{Proposition}

\newtheorem{remark}[theorem]{Remark}
\newtheorem{definition}[theorem]{Definition}

\newcommand\beq[1]{ \begin{equation}\label{#1} }
\newcommand{\eeq}{ \end{equation} }

\newcommand\beqa[1]{ \begin{eqnarray} \label{#1}}
\newcommand{\eeqa}{ \end{eqnarray} }
\newcommand{\beqano}{ \begin{eqnarray*} }
	\newcommand{\eeqano}{ \end{eqnarray*} }
\newcommand\equ[1]{{\rm (\ref{#1})}}
\newcommand{\red}{\textcolor{red}}

\def\H{{\mathcal H}}

\def\R{{\mathcal R}}

\def\real{{\mathbb R}}

\def\torus{{\mathbb T}}

\begin{document}

\title[Semi-analytical estimates for the orbital stability]
{Semi-analytical estimates for the orbital stability of Earth's satellites}

\author[I. De Blasi]{Irene De Blasi}

\address{
Department of Mathematics, University of Torino, Via Carlo Alberto 10, 00123 Torino (Italy)}
\email{irene.deblasi@unito.it}

\author[A. Celletti]{Alessandra Celletti}

\address{
Department of Mathematics, University of Roma Tor Vergata, Via
della Ricerca Scientifica 1, 00133 Roma (Italy)}

\email{celletti@mat.uniroma2.it}

\author[C. Efthymiopoulos]{Christos Efthymiopoulos}

\address{
Department of Mathematics, University of Padova, Via Trieste 63, 35121 Padova (Italy)}
\email{cefthym@math.unipd.it}

\thanks{A.C. acknowledges the MIUR Excellence
Department Project awarded to the Department of Mathematics,
University of Rome Tor Vergata, CUP E83C18000100006. A.C. and C.E.
were partially supported by EU-ITN Stardust-R and MIUR-PRIN
20178CJA2B ``New Frontiers of Celestial Mechanics: theory and
Applications''.}


\baselineskip=14pt




\begin{abstract}
Normal form stability estimates are a basic tool of Celestial Mechanics for 
characterizing the long-term stability of the orbits of natural and artificial 
bodies. Using high-order normal form constructions, we provide three different 
estimates for the orbital stability of point-mass satellites orbiting around the 
Earth. i) We demonstrate the long term stability of the semimajor axis within 
the framework of the $J_2$ problem, by a normal form construction eliminating 
the fast angle in the corresponding Hamiltonian and obtaining $\H_{J_2}$. ii) We demonstrate 
the stability of the eccentricity and inclination in a secular Hamiltonian model 
including lunisolar perturbations (the `geolunisolar' Hamiltonian $\H_{gls}$), 
after a suitable reduction of the Hamiltonian to the Laplace plane. iii) We 
numerically examine the convexity and steepness properties of the integrable 
part of the secular Hamiltonian in both the $\H_{J_2}$ and $\H_{gls}$ models, 
which reflect necessary conditions for the holding of Nekhoroshev's theorem 
on the exponential stability of the orbits. We find that the $\H_{J_2}$ model 
is non-convex, but satisfies a `three-jet' condition, while the $\H_{gls}$ 
model restores quasi-convexity by adding lunisolar terms in the Hamiltonian's 
integrable part.

\end{abstract}


\maketitle

\tableofcontents

\section{Introduction}\label{intro}
One of the major goals of Celestial Mechanics is the analysis of the stability of the 
dynamics of celestial bodies. Knowing the behavior in time of the orbital elements of an 
object allows one to predict its future, in particular whether it will cross the orbit of 
other celestial bodies and eventually undergo collisions. When applied to artificial 
spacecraft and space debris orbiting around the Earth, the question of the stability 
becomes of crucial importance, especially in view of the problem of estimating the 
orbital survival times of operating satellites or space debris. It is therefore crucial 
to devise methods that allow to study the orbital stability of objects moving around 
our planet.

In this work we will not consider the complex dynamics of an artificial spacecraft, which
should include the analysis of its shape, composition as well its rotational motion.
We will rather consider a point-mass body around the Earth, that we can identify with
one of the several millions of space debris orbiting our planet. In fact, the proliferation
of artificial satellites in the last decades has led to the generation of an
enormous amount of space debris with different sizes, from meters down to microns,
and at different altitudes. Space debris are remnants of non operational satellites or 
the result of break-up events, either collisions or explosions. Since the altitude 
determines the contribution of the different forces acting on the object (the 
gravitational attraction of the Earth, its geopotential perturbation, the influence 
of Sun and Moon, the Solar radiation pressure, etc.), it is convenient to make a 
distinction in terms of the altitude. To this end, the space in the surrounding of the 
Earth is commonly split into three main regions:
LEO (`Low Earth Orbit') denotes the region up to about 2\,000 km of altitude in which 
the Earth's attraction, the geopotential as well as the atmospheric drag are the
terms which greatly affect the dynamics of an Earth's satellite; MEO (`Medium Earth Orbit')
refers to the region between 2\,000 and 30\,000 km in which the effects of Moon, Sun and
Solar radiation pressure become important; GEO (`Geosynchronous Earth Orbit') refers to  
a thin ($\sim 200$ km) zone around the geostationary orbits (at 42\,164 km from Earth's 
center), where the satellites are in synchronous resonance with the 24-hour rotation of 
the Earth around its spin-axis.

The huge amount of objects (up to millions) in LEO, MEO, GEO needs a careful control of their orbits and 
the analysis of their dynamical stability 
(\cite{celletti2014dynamics,celletti2017dynamical,celletti2018dynamics,
CGL2020,GC2019,SARV2019}), also in view of devising appropriate mitigation measures 
(see, e.g., \cite{Jenkin,Seong,Yeon}). For objects in LEO, it is of crucial importance 
to evaluate the orbital lifetime, which is strongly affected by the atmospheric drag 
which provokes a decay of the orbits (\cite{KingHele,Krag,LinWang,Shute,Westerman}). 
In this work we focus on objects in MEO, GEO and beyond, thus not taking into account the 
dissipative effect of the atmosphere. Instead of using a propagation of the orbits 
to predict the stability time of the orbital elements, we propose a procedure based 
on analytical perturbative methods (see also \cite{gachet2017geostationary}).
More precisely, via a suitably defined sequence of canonical transformations, 
we construct a {\it normal form} of the Hamiltonian function, which enjoys the property 
that one or more of the Hamiltonian's Delaunay actions define {\it quasi-integrals 
of motion} (namely, integrals of the integrable part of the new Hamiltonian). Once the transformed Hamiltonian is obtained, the size of its {\it remainder} 
(which gives a control on the goodness of the approximation) can then be used to provide
bounds on the time variations, and hence the stability time of the orbital elements
(semimajor axis, eccentricity, inclination) as a function of the distance of the 
object from the Earth. We refer to this procedure as \textit{semi-analytical}, which means that it uses an analytical method, precisely normal forms, whose coefficients are calculated numerically,
namely with the aid of a computer.

We consider two different models to describe the motion of the debris around the Earth.
The first model takes into account only the influence of the geopotential up to the 
term $J_2$ of its expansion in spherical harmonics; we refer to this problem as the 
$J_2$ model and denote the corresponding Hamiltonian as $\mathcal{H}_{J_2}$, which results from truncating to a suitable power of the coordinates around reference values, and normalizing up to a suitable order, as described in Section \ref{ssec:nfj2}. The second model, 
referred to as the secular `geolunisolar' model (Hamiltonian $\mathcal{H}_{gls}$, truncated and normalized similarly to $\H_{J_2}$, see Sections \ref{ssec: ls finham} and \ref{ssec: NFjls}), includes also the 
effects of the Moon and the Sun, placing, for simplicity, the Moon strictly on the 
ecliptic; this last restriction means to omit from the Hamiltonian terms corresponding 
to lunisolar resonances other than the `inclination-dependent' ones. The latter 
resonances, on the other hand, are those producing the most important effects as 
regards orbital stability (see \cite{CGL2020,celletti2017dynamical} for a review). 
Furthermore, instead of formally eliminating the fast angle via canonical transformations 
(as we do in the pure $J_2$ problem), in $\mathcal{H}_{gls}$ we just take the average of the Hamiltonian 
with respect to all fast angles, namely, the mean anomaly of the satellite as well as the 
fast angles of the Moon and Sun: this averaged model allows us to focus on the satellite's long-term dynamics (i.e. the \textit{secular} one). The averaging is done in closed-form and leads to 
formulas equivalent to those described in \cite{kaula1966theory}. Furthermore, we reduce this last Hamiltonian to action-angle variables around each forced equilibrium point, which corresponds to a non-zero inclination defining the so-called \sl Laplace plane \rm (see
Section~\ref{ssec: ls finham}).

In summary, our stability estimates are obtained according to the procedure (i)-(iii) outlined below:\\
\noindent
i) Within the $J_2$ model, we make a formal elimination in the Hamiltonian of the fast 
angle (mean longitude); as a consequence, we get the preservation of the conjugate 
action variable corresponding to the semimajor axis. This allows us to compute the 
stability time for the semimajor axis at different altitudes, yielding stability times that increase with the altitude.\\ 
\noindent 
ii) Using $\mathcal{H}_{gls}$, instead, the semimajor axis becomes a parameter (with a priori 
constant value), while we proceed to analyze the behaviour of eccentricity and 
inclination. The latter is obtained using a quasi-resonant normal form, which reflects 
the 1:1 near-resonance of the integrable part of $\mathcal{H}_{gls}$ between the frequencies of the longitude of the ascending node and the sum of the argument of perigee and the longitude of the ascending node (see Section \ref{ssec: forced}). This means that, close to 
the Laplace plane, the quasi-preserved secular quantities cannot be defined neither 
as the eccentricity $e$ nor the inclination $i$ alone, but rather by the Kozai-Lidov 
combination $\mathcal{I}=1-\sqrt{1-e^2}\cos i\approx (e^2+i^2)/2$ (for $e,i$ small). We then 
explore the dependence of the stability time of $\mathcal{I}$ on the altitude of the orbit. 
Our results show that the $J_2$ and lunisolar terms have an opposite effect on the 
time of stability as the distance from the Earth increases. As a by-product of 
this analysis, we also compute the so-called forced inclination (that is, the 
inclination of the Laplace plane), which corresponds to the shift of the secular 
equilibrium from a strictly equatorial orbit to an orbit with small positive 
initial inclination, an effect caused by the fact that the perturbing bodies 
(Moon and Sun) are in orbits inclined with respect to the Earth's equator. \\
\noindent
iii) Finally, as a first step towards obtaining exponential stability estimates 
\`a la Nekhoroshev (\cite{nekhoroshev1977exponential}), we check whether some 
so-called `steepness' conditions are satisfied for the integrable part of both 
Hamiltonians $\mathcal{H}_{J_2}$ and $\mathcal{H}_{gls}$, namely whether the integrable parts are convex, 
quasi-convex, or satisfy the three-jet condition (see \cite{chierchia2018steepness} and references therein). The results show that the 
$J_2$ model is three-jet non-convex, while the contribution of the lunisolar terms 
removes the intrinsic degeneracy of the $J_2$ part and allows us to conclude that the
geolunisolar model is quasi-convex. A detailed application of the non-resonant 
form of Nekhoroshev's theorem in the Hamiltonian $\mathcal{H}_{gls}$ is the subject of 
an independent paper (see \cite{CDE2020}). 

Summarizing, the previous strategy allows us to obtain three different stability results: 
one for the semimajor axis in the $J_2$ model, a second for the stability of the 
eccentricity and inclination in the geolunisolar model, and a third on the holding, 
altogether, of necessary conditions for Nekhoroshev-type stability of the satellite 
orbits. All three results point towards the same direction, i.e. that, at least 
far from exact resonances, orbital stability can be ensured at MEO, GEO and beyond for quite long 
times ($10^4 - 10^6$ orbital periods, $10^2 - 10^4$ years). Besides these general 
numbers, one may remark that the calculation of the size of the remainder of the 
normal form actually provides an estimate of the rate of drift of the orbits in 
element space, an information  required in orbital diffusion studies for 
defunct satellites and space debris. \\

This work is organized as follows. In Section~\ref{model} we present the $J_2$ and 
geolunisolar models. Section \ref{sec:NF} briefly presents the method of the 
composition of Lie series, used in the computation of all our normal forms, 
along with some general estimates on the convergence of the normalizing transformation 
and the size of the normal form's remainder. Section~\ref{sec:stabilityJ2} focuses on 
the stability estimates with the $\mathcal{H}_{J_2}$ model, while section ~\ref{sec: stabest ls} 
deals with the stability in the framework of the $\mathcal{H}_{gls}$ model. Finally, the 
analysis of the steepness conditions for the Hamiltonians $\mathcal{H}_{J_2}$ and $\mathcal{H}_{gls}$ 
is presented in Section~\ref{sec: nondeg}. 

\section{The $J_2$ and geolunisolar models}\label{model}

Bodies orbiting around the Earth are primarily affected by the Keplerian attraction 
with our planet. However, for an accurate description of the dynamics it is mandatory 
to assume that the Earth is non-spherical. Beside the Earth, the satellite dynamics 
is subject to the gravitational influence of Sun and Moon. Section~\ref{sec:J2} 
describes the Hamiltonian model $\mathcal{H}_{J_2kep}$, which includes the Earth's Keplerian term and the  
first non-trivial term in the expansion of the geopotential. Section~\ref{ssec: ls finham} 
presents the Hamiltonian model $\mathcal{H}_{gls, sec}$, which includes $J_2$ and lunisolar terms, averaged over the fast angles. 

\subsection{The $J_2$ model}\label{sec:J2}
We consider a model describing the motion of a point-mass body, say a satellite $S$, 
under the effect of the Earth's gravitational attraction, including an approximation 
of the geopotential due to the oblateness of the Earth. Let $\underline{r}\equiv(x,y,z)$ 
be the position vector of $S$ in a geocentric reference frame, with the plane $(x,y)$ 
coinciding with the equatorial plane, and $x$ pointing towards a fixed celestial 
point (e.g. the equinox). We consider the Hamiltonian describing the motion of 
$S$ under the geopotential as the sum of two terms
\begin{equation}\label{hamJ2iniz}
\mathcal{H}_{J_2kep}=\mathcal{H}_{kep}+\mathcal{V}_{J_2}\ ,
\end{equation}
where
\begin{equation}\label{kep0}
\mathcal{H}_{kep} = {p^2\over 2} - {\mu_E\over r} 
\end{equation}
is the Keplerian part ($r=|\underline{r}|$), and 
\begin{equation}\label{hamj2fin}
\mathcal{V}_{J_2}=-J_2\frac{\mu_ER_E^2}{r^3}\left(\frac{1}{2}-\frac{3z^2}{2r^2}\right)\ ,
\end{equation}
is the $J_2$ potential term, arising from expanding the geopotential in spherical 
harmonics and retaining only the largest coefficient (see, e.g., \cite{kaula1966theory}).
The constants are the Earth's mass parameter $\mu_E = \mathcal{G}M_E$ ($\mathcal{G} = $ 
Newton's constant, $M_E = $ Earth's mass), $R_E$ is the mean Earth's radius, and 
$J_2$ is a dimensionless coefficient describing the oblate shape of the Earth. The numerical values are: 
\begin{itemize}
	\item $\mu_E=1.52984\times10^9$ $R_E^3/yr^2$; 
	\item $R_E=6378.14$ km; 
	\item $J_2= -1082.6261\times10^{-6}$.
\end{itemize}

The Hamiltonian (\ref{hamJ2iniz}) is expressed in Cartesian coordinates. However, by a standard procedure, it can be transformed 
to an expression in the following set of modified Delaunay canonical action-angle 
variables 
\begin{equation} \label{moddel}
\begin{cases}
L=\sqrt{\mu_E a}\\
P=\sqrt{\mu_E a}(1-\sqrt{1-e^2})\\
Q=\sqrt{\mu_E a}\sqrt{1-e^2}(1-\cos{i})\\
\end{cases}
\begin{cases}
\lambda=M+\omega+\Omega\\
p=-\omega-\Omega\\
q=-\Omega\ ,\\
\end{cases}
\end{equation}
where 
$(a,e,i,M,\omega,\Omega)$ are the orbital elements of the satellite (semimajor axis, 
eccentricity, inclination, mean anomaly, argument of the perigee, longitude of the 
nodes). The passage is done by first expressing the Hamiltonian (\ref{hamJ2iniz}) in 
elements via the relations (see e.g. \cite{murray1999solar})
\begin{eqnarray}\label{xyzele}
x &=&{1\over 2}r(1+\cos i)\cos(f+\omega+\Omega)
+{1\over 2}r(1-\cos i)\cos(f+\omega-\Omega)
\nonumber\\
y &=&{1\over 2}r(1+\cos i)\sin(f+\omega+\Omega)
-{1\over 2}r(1-\cos i)\sin(f+\omega-\Omega)
\\
z &=&r \sin i\sin(f+\omega),
\nonumber
\end{eqnarray}
where $f$ is the true anomaly and $r$, $\cos f$ and $\sin f$ are given by the series
\begin{eqnarray}\label{rcosfseries}
r &= &a\left[1+{e^2\over 2} - 2e\sum_{\nu=1}^{\infty}
{\left(J_{\nu-1}(\nu e)-J_{\nu+1}(\nu e)\right)\cos(\nu M)\over 2\nu} \right] 
\\ \nonumber
\cos f &= &{2(1-e^2)\over e}\sum_{\nu=1}^{\infty}J_\nu(\nu e)\cos(\nu M) - e
\\
\sin f &= &2\sqrt{1-e^2}\sum_{\nu=1}^{\infty}\frac{1}{2}
\left(J_{\nu-1}(\nu e)-J_{\nu+1}(\nu e)\right)\sin(\nu M)~~~.
\nonumber
\end{eqnarray}
In the actual calculations, all series are truncated to order $N = 15$ in the 
eccentricity $e$. Finally, we pass from the elements $(a,e,i,M,\omega,\Omega)$ 
to the canonical variables $(L,P,Q,\lambda,p,q)$ by inverting 
Eqs. (\ref{moddel}). 

To perform the high order normal form computations described in 
Section \ref{sec:stabilityJ2}, using computer algebra, it is convenient that 
the dependence of the Hamiltonian on the action-angle variables be expressed 
as a trigonometric polynomial. To this end, we first make a shift transformation 
$L\rightarrow\delta L$ around a reference value $a_*$, with 
\begin{equation}\label{deltaL}
\delta L=L-L_*=\sqrt{\mu_E a}-\sqrt{\mu_E a_*}\ .
\end{equation}
This means that the Hamiltonian found after expanding in powers of the quantity 
$\delta L$ refers to the local dynamics of orbits with semimajor axis 
$a\approx a_*$. Every time when we change the reference value $a_*$ (i.e. the 
`altitude' or `distance' of the orbit from the Earth's center), we then perform 
the Hamiltonian expansion anew around $L_*$ and obtain the stability 
estimates corresponding to that reference value. One may also note that 
$P=\mathcal{O}(e^2/2)$ and $Q=\mathcal{O}(i^2/2)$, thus all three quantities 
$\delta L$, $P$ and $Q$ are small quantities for orbits not very far from the 
equator and not very far from circular. We then expand 
$\H_{J_2kep}(\delta L,P,Q,\lambda,p,q)$ in powers of $\sqrt{\delta L}$, $\sqrt{P}$, and $\sqrt{Q}$ up 
to the same order $N=15$ as the original expansion in the eccentricity 
(this ensures missing no term in $P,Q$ in the Hamiltonian up to the order $N$). 
After this change, the truncated Hamiltonian takes the form (apart from a 
constant):
\begin{eqnarray}\label{hj2exp}
\H_{J_2}^{\leq N} &= &n_*\delta L + \omega_{1*} P +\omega_{2*} Q + \sum_{\substack{s=1\\s\neq 2}}^{2N}
\mathcal{Z}_{s}(\delta L,P,Q) \nonumber\\
~&+&\sum_{s=1}^{2N}\sum_{\substack{k_1,k_2,k_3\in\mathbb{Z}\\
0<|k_1|+|k_2|+|k_3|\leq s}}
\mathcal{P}_{s,k_1,k_2,k_3}(\delta L,P,Q)\cos(k_1\lambda+k_2p+k_3q). 
\end{eqnarray}
The functions $\mathcal{Z}_{s}$ and $\mathcal{P}_{s,k_1,k_2,k_3}(\delta L,P,Q)$ 
are polynomials of degree $s/2$ in the action variables $(\delta L,P,Q)$. Finally, the frequencies $n_*$, $\omega_{1*}$, $\omega_{2*}$ are equal 
to:
\begin{equation}\label{freqJ2}
n_* = \sqrt{{\mu_E\over a_*^3}}+J_2{3\mu_E^{1/2} R_E^2\over a_*^{7/2}},
~~~~~~~~~~~
\omega_{1*} = -J_2 {3\mu_E^{1/2} R_E^2\over 2 a_*^{7/2}},~~~~~~
\omega_{2*} = J_2 {3\mu_E^{1/2} R_E^2\over 2 a_*^{7/2}}.
\end{equation}
The Hamiltonian (\ref{hj2exp}) is the point of departure for the stability 
estimates on the orbits' semimajor axes; one notices that $\omega_{1*}= 
-\omega_{2*} = \mathcal{O}\left(J_2\right)$, a fact implying that both these 
frequencies are way smaller than $n_*\simeq (\mu_E/a_*^3)^{1/2}$ (third 
Kepler's law). Accordingly, for all orbits the angle $\lambda$ circulates 
at a rate which is $\mathcal{O}(1/J_2)$ faster than the rate of circulation 
of the angles $p,q$. Hence, $\lambda$ constitutes the `fast angle' of the 
Hamiltonian $\H_{J_2}^{\leq N}$. Its elimination through a suitable sequence of 
canonical transformations leads to the approximate constancy of the value 
of the semimajor axis, as detailed in Section \ref{sec:stabilityJ2}. 

\subsection{The geolunisolar Hamiltonian}\label{ssec: ls finham}

While stability estimates for the semimajor axis depend mostly on the Earth's 
$J_2$ term, the question of the long-term stability as regards secular variations 
in eccentricity and inclination requires considering the effects of the Lunar and
Solar gravitational tides. Let us consider a celestial body $B$ (either Moon 
or Sun) with mass $M_b$ moving around the Earth and whose orbit is exterior to that 
of the satellite. Let $\underline{r}=(x,y,z)$ and $\underline{r}_b=(x_b,y_b,z_b)$ 
be the position vectors of $S$ and $B$ in a geocentric reference frame, with 
$r=|\underline{r}|$ and $r_b=|\underline{r}_b|$. The tidal disturbance caused 
by $B$ on $S$ is described by the potential
\begin{equation}\label{bpot}
\mathcal{V}_b(\underline{r},t)=
-\mu_b\left(\frac{1}{|\underline{r}-\underline{r}_b(t)|}
-\frac{\underline{r}\cdot\underline{r}_b(t)}{r_b^3(t)}\right)
=-\frac{\mu_b}{r_b(t)}+\frac{\mu_b}{2r_b^3(t)}r^2+\frac{3}{2}\frac{\mu_b(\underline{r}\cdot\underline{r}_b(t))^2}{r_b^5(t)}+\mathcal{O}\left(\left(\frac{r}{r_b}\right)^3\right)\ ,
\end{equation}
where $\mu_b=\mathcal{G}M_b$.
The first term $-\mu_b/r_b$ in the multipolar expansion (\ref{bpot}) does not 
depend on the coordinates of $S$, therefore it can be omitted from the Hamiltonian 
of motion of $S$. Thus, the tidal (or `third body') perturbation terms in the 
Hamiltonian takes the form: 
\begin{equation}\label{htide}
\H_{3B} = 
\mu_m\left({1\over 2}\frac{r^2}{r_m^3(t)}
-{3\over 2}\frac{(\underline{r}\cdot \underline{r}_m(t))^2}{r_m^5(t)}\right)
+
\mu_s\left({1\over 2}\frac{r^2}{r_s^3(t)}
-{3\over 2}\frac{(\underline{r}\cdot \underline{r}_s(t))^2}{r_s^5(t)}\right)
+O_3=\H_m+\H_s, 
\end{equation}
where $\mu_m,\underline{r}_m$ and $\mu_s,\underline{r}_s$ are the mass and geocentric 
position vectors of the Moon and Sun respectively, and $O_3$ denotes octupolar or 
higher order terms in the expansion of the third body potentials. The exact form 
of the term $\H_{3B}$ depends now on the model adopted for the geocentric orbits of 
the Sun and Moon. In the framework of the present paper, we adopt the following 
models for Sun and Moon: 
\begin{enumerate}
    \item we suppose that the Sun's geocentric orbit is an ellipse 
    lying in the Earth's ecliptic plane
    (i.e., with inclination $i_{s0}=23.43^\circ$ with respect to the equatorial plane),
    $\Omega_s=0^\circ$, $a_s=1.496\cdot 10^8\, km$ and $e_s=0.0167$;
    \item we assume the Lunar orbit as elliptic and also lying on the ecliptic plane, 
    with $a_m=384748$ $km $, $e_m=0.065$ and $i_{m0}=i_{s0}$. Note that this assumption 
    ignores the precession of the Lunar node (with period $\simeq 18.6$ yr) associated 
    with the inclination of the Moon's orbit with respect to the ecliptic 
    (by $5^\circ15'$). While the precession of the Lunar node is important near 
    secular lunisolar resonances\footnote{By secular lunisolar resonances we mean resonances of the form $k_1\dot{\omega}+k_2\dot{\Omega}+k_3\dot{\Omega}_M=0$, with $(k_1,k_2,k_3)\in\mathbb{Z}^3\backslash\{\underline{0}\}$, thus involving the rate of variation of the longitude of the ascending node of the Moon.}, it only has a minimal effect far from these resonances, 
    as substantiated by numerical studies (e.g. \cite{gkolias2016order}, \cite{rosengren2018dynamical}). Thus, we ignore this effect 
    in our present estimates (Section 6).
\end{enumerate}

Under the above approximations, the satellite Hamiltonian $\mathcal{H}_{J_2ls}$ takes the form 
\begin{equation}\label{hamJ2ls}
\mathcal{H}_{J_2ls}=\mathcal{H}_{J_2kep}+\H_{3B}.\
\end{equation}
This is a Hamiltonian depending on three degrees of freedom (the coordinates and 
momenta of the satellite) as well as on time (through the vectors $\underline{r}_m(t)$ 
and $\underline{r}_s(t)$). However, contrary to the case of the Hamiltonian $\H_{J_2kep}$, 
in which we are interested in establishing the long-term stability of the semimajor 
axis over short-period oscillations, here we are interested in the question of the
stability of the eccentricity and inclination of the satellite over secular timescales. 
Thus, as customary (see \cite{kaula1966theory}, \cite{celletti2017analytical}), we average $\H_{J_2ls}$ with respect to the mean anomalies 
of the satellite, Moon and Sun. The averaging can be done in closed form (see, for 
example, \cite{kaula1966theory}), and leads to:
$$
\H^{(av)}_{J_2}={1\over 2\pi}\int_0^{2\pi}\H_{J_2kep} dM = 
\int_0^{2\pi}\H_{J_2kep}{r^2\over a^2\sqrt{1-e^2}}df
$$
$$
\H^{(av)}_{m}={1\over 4\pi^2}\int_0^{2\pi}\int_0^{2\pi}\mathcal{V}_{m} dM dM_m = 
\int_0^{2\pi}\int_0^{2\pi}
\mathcal{V}_{m}(1-e\cos E){r_m^2\over a_m^2\sqrt{1-e_m^2}}dEdf_m
$$
$$
\H^{(av)}_{s}={1\over 4\pi^2}\int_0^{2\pi}\int_0^{2\pi}\mathcal{V}_{s} dM dM_s = 
\int_0^{2\pi}\int_0^{2\pi}
\mathcal{V}_{s}(1-e\cos E){r_s^2\over a_s^2\sqrt{1-e_s^2}}dEdf_s.
$$
Here, $f,E$ are the satellite's true and eccentric anomaly, while $f_m,f_s$ are 
the Moon's and Sun's true anomaly along their geocentric orbits. The averaged $J_2$ 
term takes the form (apart from constant terms):
\begin{equation}\label{J2}
\mathcal{H}_{J_2}^{(av)}(e,i,\omega,\Omega)=-J_2\frac{\mu_E
R_E^2}{a^3(1-e^2)^{3/2}}\left(\frac{1}{2}-\frac{3}{4}\sin^2{i}\right).
\end{equation}
The terms $\H^{(av)}_{m}(e,i,\omega,\Omega)$ and $\H^{(av)}_{s}(e,i,\omega,\Omega)$, 
instead, turn out to be identical to those given in equations (3.6) and (3.7) of 
\cite{celletti2017dynamical}, setting $i_M=0$. Then, the Hamiltonian averaged over all short 
period terms, hereafter referred to as the \textit{secular geolunisolar 
Hamiltonian}, takes the form
\begin{equation}\label{hamavfin}
\mathcal{H}_{gls, sec}(e,i,\omega,\Omega)=\mathcal{H}_{J_2}^{(av)}
+\mathcal{H}_{s}^{(av)}+\mathcal{H}_{m}^{(av)}\ ,
\end{equation}
which, in terms of the Delaunay modified variables, has two degrees of freedom. 

\subsubsection{Expansion around the forced inclination} \label{ssec: forced}
As it was done in the case of the $J_2$ model (Eq. (\ref{hj2exp})), normal form 
computations for the Hamiltonian (\ref{hamavfin}), expressed in Delaunay action-angle 
variables, require a polynomial expansion in the action variables around some 
preselected values. In the case of the secular geolunisolar Hamiltonian (\ref{hamavfin}), 
a natural choice of the origin for such expansions is the {\it forced element} values: 
writing $\H_{gls, sec}$ as a function of the Delaunay variables, say, $\H_{gls, sec}(P,Q,p,q;a)$ 
(where the semimajor axis $a$ is now a priori constant, hence, can be considered 
as a parameter), a forced equilibrium is defined as an equilibrium point of the 
secular Hamiltonian, i.e., a point $(Q^{(eq)},P^{(eq)},q^{(eq)},p^{(eq)})$ for 
which the following relations hold:
\begin{equation}\label{eqforced}
\left({\partial \H_{gls, sec}\over\partial P}\right)_{eq}
=
\left({\partial \H_{gls, sec}\over\partial Q}\right)_{eq}
=
\left({\partial \H_{gls, sec}\over\partial p}\right)_{eq}
=
\left({\partial \H_{gls, sec}\over\partial q}\right)_{eq}
=0,
\end{equation}
where the subscript `eq' denotes the condition 
$Q=Q^{(eq)},P=P^{(eq)},q=q^{(eq)},p=p^{(eq)}$. 
In the case of the Hamiltonian (\ref{hamavfin}), a forced equilibrium solution can 
be computed by writing first $\H_{gls, sec}$ in terms of the \textit{Poincar\'e 
variables} as
\begin{equation}\label{possvar}
\begin{cases}
X_1=\sqrt{2Q}\sin{q}\ , \quad &X_2=\sqrt{2P}\sin{p}\ ,\\
Y_1=\sqrt{2Q}\cos{q}\ , \quad &Y_2=\sqrt{2P}\cos{p}\ .
\end{cases}
\end{equation}
Expanding up to quadratic terms in the Poincar\'{e} variables, the truncated secular 
Hamiltonian has the form
\beq{Htilde}
\widetilde{\mathcal{H}}(Y_1,Y_2,X_1,X_2)=A_1Y_1+B_1(X_1^2+Y_1^2)+B_2(X_2^2+Y_2^2)\ ,
\eeq
where the coefficients $A_1$, $B_1$ are given by:
\begin{equation}\label{A1B1}
\begin{split}
&A_1=-\frac{3R_E^2a^{7/4}\sin{(2i_0)}}{8(\mathcal{G}M_E)^{1/4}}
\left(\frac{\mathcal{G}M_m}{a_m^3}+\frac{\mathcal{G}M_s}{a_s^3}\right),\\
&B_1=\frac{3}{4}\frac{\sqrt{\mathcal{G}M_E}R_E^2J_2}{a^{7/2}}
+\frac{3\mathcal{G}M_m(2-3\sin^2{i_0})}{16\sqrt{\frac{\mathcal{G}M_E}{a^3}}a_m^3}
+\frac{3\mathcal{G}M_s(2-3\sin^2{i_0})}{16\sqrt{\frac{\mathcal{G}M_E}{a^3}}a_s^3}.
\end{split}
\end{equation}
The Hamiltonian (\ref{Htilde}) corresponds to two decoupled harmonic oscillators in 
the variables $(X_1, Y_1)$ and $(X_2, Y_2)$. The second harmonic oscillator 
(corresponding to the action-angle pair $(P,p)$, hence, to the orbit's eccentricity 
vector) has an equilibrium point at $\left(X_2^{(eq)},Y_2^{(eq)}\right)=(0,0)$, 
implying $P^{(eq)}=0$ and any value $0\leq p^{(eq)}<2\pi$. 
This implies that the sub-manifold of circular orbits $e=0$ (corresponding to $P=0$) 
is invariant under the flow of the secular geolunisolar Hamiltonian. On the other 
hand, as regards the pair $(X_1,Y_1)$, Hamilton's equations for the Hamiltonian  
\equ{Htilde} yield:
\begin{equation}\label{hameqosc}
\begin{cases}
\dot{X_1}=A_1+2B_1Y_1 \\
\dot{Y_1}=-2B_1X_1\ .
\end{cases}
\end{equation}
For $i_0\neq0$, the equilibrium point of (\ref{hameqosc}) is given by
\begin{equation*}
X_1^{(eq)}=0\ , \quad Y_1^{(eq)}=-\frac{A_1}{2B_1}\neq 0.
\end{equation*}
Setting 
$Q^{(eq)} = \left(\left(X_1^{(eq)}\right)^2+\left(Y_1^{(eq)}\right)^2\right)/2$, 
$i_{eq} \simeq (2Q^{(eq)}/\sqrt{\mu_E a})^{1/2}$ (for $Q^{(eq)}$ small), we arrive 
at
\begin{equation}\label{finc}
i^{(eq)}\simeq-\frac{A_1}{2B_1}\frac{1}{(\mu_E a)^{1/4}},
~~~q^{(eq)}=-\Omega^{(eq)} = 0.
\end{equation}
More accurate expressions for the {\it forced inclination} $i^{(eq)}$ can be 
obtained by introducing (\ref{finc}) along with the remaining equilibrium values 
in the derivatives of the full secular Hamiltonian (\ref{hamavfin}) and finding 
the roots of Hamilton's equations. One can readily verify that $q^{(eq)}=0$ at 
all orders, while $i^{(eq)}$ is subject to small corrections with respect to 
the expression (\ref{finc}). In physical terms, the forced inclination $i^{(eq)}$ 
defines the inclination of the {\it Laplace plane}: since the perturbing bodies 
(Moon and Sun) are in orbits inclined with respect to the equator, a satellite 
orbit can maintain its inclination constant when the latter has the value 
$i^{(eq)}$. Inspecting the form of the coefficients (\ref{A1B1}), we find that 
$i^{(eq)}\rightarrow 0$ as $a\rightarrow 0$, while it can be shown that $i^{(eq)}\rightarrow i_0$ for values of $a$ greater than the GEO one (see for example \cite{rosengren2014classical}), reflecting the fact that the Laplace plane tends to 
coincide with the equator for satellite orbits close to the Earth (as imposed 
by the oblateness of the Earth), while it tends to coincide with the ecliptic 
at large distances from the Earth (where the Lunar and Solar tides dominate). 

Returning to the expansion of the secular geolunisolar Hamiltonian, making the 
shift transformation $\delta Y_1=Y_1-Y_1^{(eq)}$ allows us to express the 
Hamiltonian as a polynomial in the variables $(X_1,\delta Y_1)$ and $(X_2,Y_2)$. 
The Hamiltonian $\H_{gls, sec}$ starts now with terms of second degree which we regroup in $\H_2$: 
\begin{equation}\label{hamquad}
\mathcal{H}_2= 
{b_1+\epsilon_1\over 2} X_1^2 + {b_1+\epsilon_2\over 2} \delta Y_1^2 +
{b_1+\epsilon_3\over 2} X_2^2 +{b_1+\epsilon_4\over 2} Y_2^2,
\end{equation}
where $b_1 = 2B_1$ and $\epsilon_1$, $\epsilon_2$, $\epsilon_3$, $\epsilon_4$ are corrections of order 
$\mathcal{O}(\mu_b a^{3/2}/(\mu_E^{1/2}a_b^3))$, with the index $b$ referring to 
the Moon or Sun. All these corrections turn to be rather small, with relative 
size $\sim 10^{-3} B_{10}$ at semimajor axis $a\sim 10^4$ km, where
$$
B_{10} = \frac{3}{4}\frac{\sqrt{\mathcal{G}M_E}R_E^2J_2}{a^{7/2}}.
$$
Thus, after a canonical rescaling $X_1 = c_{12} \tilde{X}_1$, $\delta Y_1=\tilde{Y}_1/c_{12}$, $X_2 = c_{34} \tilde{X}_2$, $ Y_2=\tilde{Y}_2/c_{34}$,
with $(c_{12})^4 = (b_1+\epsilon_2)/(b_1+\epsilon_1) = 
1 + \mathcal{O}(\mu_b a^{3/2}/(B_{10}\mu_E^{1/2}a_b^3)$ and $(c_{34})^4 = (b_1+\epsilon_4)/(b_1+\epsilon_3) = 
1 + \mathcal{O}(\mu_b a^{3/2}/(B_{10}\mu_E^{1/2}a_b^3)$, the secular lunisolar 
Hamiltonian resumes the form:
\begin{eqnarray}\label{hamavfinpol}
\H_{gls, sec} 
&= &{\nu_1\over 2}\left(\tilde{X}_1^2+\tilde{Y}_1^2\right) 
+ {\nu_2\over 2}\left(\tilde{X}_2^2+\tilde{Y}_2^2\right) \\
&+& \sum_{s=3}^\infty
\sum_{\substack{k_1,k_2,l_1,l_2\in\mathbb{N}\\
k_1+k_2+l_1+l_2=s}}
h_{k_1,k_2,l_1,l_2}\tilde{X}_1^{k_1}\tilde{X}_2^{k_2}\tilde{Y}_1^{l_1}\tilde{Y}_2^{l_2}.
\nonumber
\end{eqnarray}
This is the typical form of a secular Hamiltonian, consisting of linear oscillators 
(with frequencies $\nu_1$, $\nu_2$) coupled with nonlinear terms. However, we 
have $\nu_1=\nu_2 + \mathcal{O}(\mu_b a^{3/2}/\mu_E^{1/2}a_b^3)$, implying that 
the two frequencies are nearly equal 
$$
\nu_1\simeq \nu_2 \simeq \frac{3}{2}\frac{\sqrt{\mathcal{G}M_E}R_E^2J_2}{a^{7/2}}~~.
$$
This is a consequence of the axisymmetry of the $J_2$ model, implying that the 
secular frequencies $\dot{q}=-\dot{\Omega}$ and $\dot{p}=-\dot{\omega}-\dot{\Omega}$ 
are equal for nearly equatorial orbits in this model. As we will see in Section 
\ref{sec: stabest ls}, this near-equality implies that with the present normal 
form estimates one cannot establish independently the long term stability of 
the eccentricity and the inclination, but only the long-term stability of 
the Kozai-Lidov integral $\mathcal{I} = \tilde{X}_1^2+\tilde{Y}_1^2 + \tilde{X}_2^2 + \tilde{Y}_2^2$, which 
couples oscillations between the eccentricity and the proper inclination of the 
satellite. 

Finally, the Hamiltonian (\ref{hamavfinpol}) can be written in action-angle 
variables 
$\tilde{X}_1 = \sqrt{2I_1}\sin\phi_1$, $\tilde{Y}_1 = \sqrt{2I_1}\cos\phi_1$, 
$\tilde{X}_2 = \sqrt{2I_2}\sin\phi_2$, $\tilde{Y}_2 = \sqrt{2I_2}\cos\phi_2$ as
\begin{eqnarray}\label{hamavfinaa}
\H_{gls, sec} 
&= &\nu_1 I_1 + \nu_2 I_2 \\
&+& \sum_{s=3}^\infty
\sum_{\substack{s_1,s_2\in\mathbb{N}\\
s_1+s_2=s}}
\sum_{\substack{k_1,k_2\in\mathbb{Z}\\
|k_1|+|k_2|\leq s\\
(|k|_1+|k|_2)\text{ }\equiv\text{ } s \text{ }(\text{mod }2)}}
\tilde{h}_{s_1,s_2,k_1,k_2} I_1^{s_1/2}I_2^{s_2/2}\cos(k_1\phi_1+k_2\phi_2). 
\nonumber
\end{eqnarray}
The Hamiltonian (\ref{hamavfinaa}) is the starting point for all normal form 
calculations in Section \ref{sec: stabest ls}. For computational reasons, the expansion in (\ref{hamavfinaa}) is truncated up to a maximal order $N=15$, leading to the truncated form 
\begin{eqnarray}\label{hamavfinaaN}
\H_{gls, sec}^{\leq N}(I_1,I_2,\phi_1,\phi_2) 
&= &\nu_1 I_1 + \nu_2 I_2 \\
&+& \sum_{s=3}^N
\sum_{\substack{s_1,s_2\in\mathbb{N}\\
		s_1+s_2=s}}
\sum_{\substack{k_1,k_2\in\mathbb{Z}\\
		|k_1|+|k_2|\leq s\\
		(|k|_1+|k|_2)\text{ }\equiv\text{ } s \text{ }(\text{mod }2)}}
\tilde{h}_{s_1,s_2,k_1,k_2} I_1^{s_1/2}I_2^{s_2/2}\cos(k_1\phi_1+k_2\phi_2). 
\nonumber
\end{eqnarray}
\section{Hamiltonian Normalization}\label{sec:NF}

In this Section we briefly recall some basic definitions related to normal form theory 
and its use in obtaining stability estimates based on the size of the normal form's 
\textit{remainder}. In Sections \ref{sec:stabilityJ2} and \ref{sec: stabest ls} we 
will discuss the particular normalizations implemented on the Hamiltonians 
(\ref{hj2exp}) and (\ref{hamavfinaaN}) respectively.  

\subsection{Normal form and remainder}
Consider a Hamiltonian function of the form 
\begin{equation}\label{ham0}
\H(\underline{A},\underline{\varphi})
= Z_0(\underline{A})
+ \H_1(\underline{A},\underline{\varphi})
=\underline{\omega}\cdot\underline{A}
+ \H_1(\underline{A},\underline{\varphi}),
\end{equation}
where $\omega_j$ are real constants, and 
$(\underline{A},\underline{\varphi})\in\real^n\times\torus^n$ are action-angle 
variables. We assume that $\H_1$ is analytic in the complex domain 
$(\underline{A},\underline{\varphi})\in$ 
$D_{\rho,\sigma}(U)=B_\rho U\times S_\sigma$ 
(or simply $D_{\rho,\sigma}$), where $U$ is an open domain of $\real^n$, 
$B_\rho U$ is a complex neighborhood of $U$ of size $\rho$:
\begin{equation}\label{complexaction}
B_\rho U=\{\underline{A}\in\mathbb{C}^n:\ dist(\underline{A},U)<\rho\}\ ,
\end{equation}
$S_\sigma$ is the complex strip
\begin{equation}\label{complexangles}
S_\sigma=\{\underline{\varphi}\in\mathbb{C}^n:\ Re(\varphi_j)\in\mathbb{T},\
\vert Im(\varphi_j)\vert<\sigma, \quad j=1,\dots,n\}
\end{equation} 
for $\rho,\sigma>0$.
On $D_{\rho,\sigma}(U)$ we define the norm of a function 
$f=f(\underline{A},\underline{\varphi})$ as
\begin{equation}\label{normnf}
\|f\|_{\rho,\sigma}=\sup_{(\underline{A},\underline{\varphi})\in 
D_{\rho,\sigma}}|f(\underline{A},\underline{\varphi})|\ .
\end{equation}
The aim of normalization theory  is to introduce a near to identity canonical 
transformation 
$\Phi:~
(\underline{A},\underline{\varphi})\rightarrow
(\underline{A}',\underline{\varphi}')$, 
so that in the new variables $(\underline{A}',\underline{\varphi}')$ the 
Hamiltonian (\ref{ham0}) takes the form
\begin{equation}\label{hamnf}
\H\left(\underline{A}(\underline{A}',\underline{\varphi}'),
\underline{\varphi}(\underline{A}',\underline{\varphi}')\right) = 
Z(\underline{A}',\underline{\varphi}')+R(\underline{A}',\underline{\varphi}')
\end{equation}
with the following properties:\\
i) the transformation $\Phi$ is analytic in a domain $D_{\rho',\sigma'}(U)$ 
with $0<\rho'<\rho$, $0<\sigma'<\sigma$,\\
ii) the dynamics under $Z(\underline{A}',\underline{\varphi}')$, called the 
\textit{normal form}, has some desired properties (see below), and \\
iii) under the norm definition (\ref{normnf}) one has $\|R\|_{\rho',\sigma'}\ll
\|Z\|_{\rho',\sigma'}$ implying that the function 
$R(\underline{A}',\underline{\varphi}')$, called the \textit{remainder}, 
introduces only a small correction with respect to the flow under the normal 
form term 
$Z(\underline{A}',\underline{\varphi}')$.

Regarding point ii) above, see, e.g., \cite{efthymiopoulos2011canonical} for a
definition of the properties of the normal form term in various contexts of 
perturbation theory (e.g. in the Kolmogorov-Arnold-Moser or Nekhoroshev theories). 
Here we mention three cases of particular interest, pertinent to our present work:\\

\noindent
\textit{Case 1: Birkhoff normal form}. The function $Z$ is independent of the 
angles $\underline{\varphi}'$. This kind of normalization allows us to prove the 
near-constancy of the action variables $\underline{A}$.\\
\\
\textit{Case 2: elimination of short-period terms}. The real constants $\omega_j$ 
in (\ref{ham0}), called the unperturbed frequencies, are divided in two groups,  
`fast' $\{\omega_1,\ldots,\omega_{K_f}\}$, $1\leq K_f<n$, and `slow' 
$\{\omega_{K_f+1},\ldots,\omega_{n}\}$, such that 
$\min\{|\omega_1|,\ldots,|\omega_{K_f}|\}\gg
\max\{|\omega_{K_f+1}|,\ldots,|\omega_{n}|\}$.
In this case, it turns convenient to introduce a normalizing transformation $\Phi$ 
such that the normal form $Z$ becomes independent of the `fast angles' 
$\{\varphi_1',\ldots,\varphi_{K_f}'\}$. Such is the case of the normal form 
encountered in Section \ref{sec:stabilityJ2}, leading to estimates on the stability 
of the semimajor axis in the $J_2$ problem. The corresponding Hamiltonian is of the 
form (\ref{ham0}), with $n=3$, $A_1=\delta L$, $A_2 = P$, $A_3=Q$, $\varphi_1=\lambda$, 
$\varphi_2=p$, $\varphi_3=q$.\\
\\
\textit{Case 3: resonant normal form}. The frequencies $\omega_j$ satisfy one or 
more quasi-commensurability relations of the form $\underline{m}\cdot\underline{\omega}
\simeq 0$, with $\underline{m}\in\mathbb{Z}^n,|\underline{m}|\neq 0$. The maximum number 
of linearly independent and irreducible integer vectors $\underline{m}_l$, $1\leq l\leq 
l_{max}$, yielding exact commensurabilities for a given set of frequencies $\omega_j$, 
satisfies $0\leq l_{max}\leq n$. Since $\H_1$ is analytic in $D_{\rho,\sigma}(U)$ 
and periodic in $\underline{\varphi}$, $\H_1$ admits the Fourier decomposition
\begin{equation}\label{h1four}
\H_1(\underline{A},\underline{\varphi}) = 
\sum_{\underline{k}\in\mathbb{Z}^n} h_{1,\underline{k}}(\underline{A})
e^{i\underline{k}\cdot\underline{\varphi}},
\end{equation}
where, according to Fourier theorem, the coefficients 
$|h_{1,\underline{k}}(\underline{A})|$ are bounded by exponentially decaying 
quantities $\mathcal{O}(e^{-|\underline{k}|\sigma})$. Then, it turns out that the 
appropriate normal form $Z$ has the \textit{resonant form}:
\begin{equation}\label{nfres}
Z(\underline{A}',\underline{\varphi}') = 
\sum_{\underline{k}\in\mathcal{M}} 
\zeta_{\underline{k}}(\underline{A'})
e^{i\underline{k}\cdot\underline{\varphi}'},
\end{equation}
for some Fourier coefficients $\zeta_{\underline{k}}(\underline{A}')$ and where
$$
\mathcal{M}:=\{\underline{k}\in\mathbb{Z}: \underline{k}\cdot\underline{m}_l=0 
\mbox{~for all~} l=1,\ldots,l_{max}\}
$$
is the `resonant module'. A normal form of the form (\ref{nfres}) implies the 
existence of $n-l_{max}$ quasi-integrals of the form $I_i = \underline{K}_i\cdot
\underline{A}$, $i=1,\ldots,n-l_{max}$, where the vectors $\underline{K}_i$ 
satisfy the equations $\underline{K}_i\cdot\underline{m}_l=0$ for all $l$ with 
$1\leq l\leq l_{max}$. The quantities $I_i$ are called the \textit{resonant integrals}
of the Hamiltonian (\ref{nfres}). 

As an example, whenever $\nu_1=\nu_2$ the secular geolunisolar Hamiltonian (\ref{hamavfinaaN}) admits a 
resonant normal form. We have $n=2$, $l_{max}=1$, $\underline{m}_1 = (1,-1)$, 
$A_1=I_1$, $A_2=I_2$, $\varphi_1=\phi_1$, $\varphi_2 = \phi_2$. Therefore, the 
normal form contains terms independent of the angles or depending on the angles 
through trigonometric terms of the form $\cos(k(\phi_1'-\phi_2'))$, $k=1,2,\ldots$. 
The associated resonant integral corresponds to the `Kozai-Lidov' integral 
$\mathcal{I}=I_1+I_2$ (see \cite{Lidov}).

\begin{definition}\label{nfconv}

A $r$-th step Hamiltonian normalization process is a composition of near identity 
transformations
\begin{equation}\label{rseq}
\Phi^{(r)} = \Phi_r\circ\Phi_{r-1}\circ\ldots\circ\Phi_1
\end{equation}
mapping the initial action-angle variables to the r-th step normalized action-angle 
variables via the successive transformations $(\underline{A}^{(s)},\underline{\phi}^{(s)})=\Phi^{(s)}(\underline{A}^{(s-1)},\underline{\phi}^{(s-1)})$, 
$s=1,2,,\ldots,r$, $(\underline{A}^{(0)},
\underline{\varphi}^{(0)})\equiv(\underline{A},\underline{\varphi})$, defined so 
that the compositions
$$
\Phi^{(s)} = \Phi_s\circ\Phi_{s-1}\circ\ldots\circ\Phi_1
$$
for all $s=1,\ldots,r$ are analytic and with inverse analytic within non-null domains $D_{\rho^{(s)},\sigma^{(s)}}
\neq\emptyset$, and the rth-step Hamiltonian takes the form:
\begin{equation}\label{nfrstep}
\H^{(r)}\left(\underline{A}(\underline{A}^{(r)},\underline{\varphi}^{(r)}),
\underline{\varphi}(\underline{A}^{(r)},\underline{\varphi}^{(r)})\right) = 
Z^{(r)}(\underline{A}^{(r)},\underline{\varphi}^{(r)})+R^{(r)}(\underline{A}^{(r)}
,\underline{\varphi}^{(r)})
\end{equation}
with $\|R^{(r)}\|_{\rho^{(r)},\sigma^{(r)}}\ll\|Z^{(r)}\|_{\rho^{(r)},\sigma^{(r)}}$.
\end{definition}

The semi-analytical estimates of stability that we will develop in the next sections are based on defining a suitable $r$-step 
sequence of canonical transformations $\Phi_1,\Phi_2,\ldots,\Phi_r$ reducing the size of 
the remainder $\|R^{(r)}\|_{\rho^{(r)},\sigma^{(r)}}$ as much as possible given the 
initial Hamiltonian model considered. The appropriate sequence is found using the method 
of \textit{Lie series} (see Section \ref{sec:stabilityJ2}). The obtained times of 
stability are of order $\|R^{(r_{opt})}\|_{\rho^{(ropt)},\sigma^{(ropt)}}^{-1}$, 
where $r_{opt}$ is the normalization order yielding the smallest possible remainder norm.  
The value of $r_{opt}$ can be obtained via theoretical estimates (see \cite{efthymiopoulos2004nonconvergence}), but in 
practice, it is also limited by the maximum order in which our computer-algebra 
normal form calculations can proceed. Theoretical estimates imply that the size of the 
remainder norm is exponentially small in the inverse of the size of the perturbation 
$\|\H_1\|_{\rho,\sigma}$ in Eq. (\ref{ham0}). For example, in the simplest case of the Birkhoff normal 
form, we have the following theorem (see \cite{fasso1989composition} for full details).

\begin{theorem}\label{fassothm}
Consider the Hamiltonian expressed in action-angle variables
$\mathcal{H}(\underline{A},\underline{\varphi})=
\underline{\omega}\cdot\underline{A}+f(\underline{A},\underline{\varphi})$,
where $\underline{\omega}\in\mathbb{R}^n$ satisfies the following Diophantine 
condition: there exist $\tau,\gamma>0$ such that
\begin{equation}\label{DC}
|\underline{k}\cdot\underline{\omega}|\geq\frac{\gamma}{|\underline{k}|^\tau}\
\qquad
\forall\underline{k}\in\mathbb{Z}^n\backslash{\{\underline{0}\}}\
\end{equation}
and $f$ is real analytic on $D_{\rho,\sigma}$ for some $\rho,\sigma>0$.
Consider two positive parameters $\delta<\rho/2$ and $\xi<\sigma/2$, and for
$r\geq1$, let
\begin{equation}\label{fassoeq1}
\epsilon_1^*=\frac{\gamma\delta\xi^{\tau+1}}{2^{n-\tau+4}\sqrt{(2\tau+2)!}\
\|f\|_{\rho,\sigma}}\ ,\qquad
\epsilon_r^*=\frac{\epsilon_1^*}{r^{\tau+2}}\ .
\end{equation}
Then, for any
\begin{equation}\label{fassoeq3}
r<\left(\frac{\gamma\delta\xi^{\tau+1}}{2^{n-\tau+4}
\sqrt{(2\tau+2)!}}\right)^{1/\tau+2}\frac{1}
{\| f\|_{\rho,\sigma}^{1/\tau+2}}\ ,
\end{equation}
there exists a real analytic canonical transformation
$\Phi:D_{\rho-2\delta, \sigma-2\xi}\mapsto D_{\rho,\sigma}$ such
that the transformed Hamiltonian has the form
\begin{equation}
\mathcal{H}\circ\Phi=h(\underline{A})+
\sum_{s=1}^rZ_s(\underline{A})
+\red{R}^{(r+1)}(\underline{A},\underline{\varphi})\
,
\end{equation}
where the remainder $R^{(r+1)}$ can be bounded as
\begin{equation}\label{fassoeq5}
\| R^{(r+1)}\|_{\rho-2\delta, \sigma-2\xi}\leq
\frac{\|f\|_{\rho,\sigma}}{4r^{\tau+2}}
\left(\frac{1}{\epsilon_r^*}\right)^r\frac{\epsilon_r^*}{\epsilon_r^*-1}\ .
    \end{equation}
\end{theorem}

Casting together (\ref{fassoeq1}) and (\ref{fassoeq5}), one readily sees that the 
remainder grows more rapidly than any power of $r$, namely as $(r^{\tau+2})^{r-1}$. 
Consequently, this procedure does not converge for $r\rightarrow\infty$. In any case, 
we remark that, as the threshold value for the normalization order $r$ is proportional 
to the inverse of $\|f\|_{\rho,\sigma}^{1/(\tau+2)}$, if we manage to reduce the size 
of the initial remainder function, then we can increase the maximum value of $r$ for 
which Theorem \ref{fassothm} is satisfied.

Similar estimates hold in the case of the resonant normal form constructions (see \cite{efthymiopoulos2004nonconvergence}).  
The behavior of the size of the remainder as a function of the normalization order 
$r$ will be examined in detail in our semi-analytical computations in Sections 
\ref{sec:stabilityJ2} and \ref{sec: stabest ls} below. 

\subsection{Book-keeping and construction of the normal form}
\label{ssec: j2 gen}
Both Hamiltonians (\ref{hj2exp}) and (\ref{hamavfinaaN}) are of the form (\ref{ham0}), 
therefore the above results on Hamiltonian normalization apply. In order to compute 
the composition of canonical transformations required in Eq. (\ref{rseq}), we 
implement the method of \textit{composition of Lie series}, after introducing 
a suitable {\textit{book-keeping} (see \cite{efthymiopoulos2011canonical}) to 
separate terms in the Hamiltonian according to estimates of their order of 
smallness.  

\begin{definition}
Consider a `book-keeping symbol' $\epsilon$, with numerical value $\epsilon=1$. 
A book-keeping rule is a splitting of the initial Hamiltonian $\H(\underline{A},
\underline{\varphi})$ in the form
\begin{equation}\label{bkeep}
\H(\underline{A},\underline{\varphi}) = \underline{\omega}\cdot\underline{A}
+ \sum_{s=1}^\infty \epsilon^s \H_s(\underline{A},\underline{\varphi}).
\end{equation}
\end{definition}

\begin{remark}
The splitting can in principle be arbitrary. However, the sequence of remainders 
$\|R^{(r)}\|_{\rho^{(r)},\sigma^{(r)}}$ found by Hamiltonian normalization behaves 
well, i.e. $\|R^{(s)}\|_{\rho^{(s)},\sigma^{(s)}}<\|R^{(s-1)}\|_{\rho^{(s-1)},
\sigma^{(s-1)}}$ for $s=1,\ldots,r_{opt}$ when the splitting (\ref{bkeep}) is done 
so as to reflect the order of smallness of different terms in the 
Hamiltonian. Roughly speaking, one must have $\|\H_s\|_{\rho,\sigma} = 
\mathcal{O}\left(\|\H_1\|_{\rho,\sigma}^s\right)$
(see \cite{efthymiopoulos2011canonical}).
\end{remark}

\begin{proposition}\label{lieser}
{\bf Lie series:} Let $\chi(\underline{A},\underline{\varphi})$, called the Lie 
generating function, be a function analytic in the domain $D_{\rho,\sigma}(U)$, 
and $\mathcal{L}_\chi$ denote the Poisson bracket operator $\mathcal{L}_\chi\cdot 
= \{\cdot,\chi\}$. Given positive numbers $\delta<\rho$ and $\xi<\sigma$, assume 
that 
$$
\min\left(
\delta\|{\partial\chi\over\partial q}\|_{\rho-\delta,\sigma-\xi}^{-1},
\xi\|{\partial\chi\over\partial p}\|_{\rho-\delta,\sigma-\xi}^{-1}\right)>1.~~
$$
Then, the mapping 
\begin{equation}\label{lietra}
(\underline{A}',\underline{\varphi}') = 
\exp(\mathcal{L}_\chi)(\underline{A},\underline{\varphi}) = 
\sum_{j=0}^\infty {1\over j!}\mathcal{L}_\chi^j(\underline{A},\underline{\varphi})
\end{equation}
is an analytic canonical transformation of the domain $D_{\rho-\delta,\sigma-\xi}(U)$ 
onto itself.
\end{proposition}

The proof consists in implementing Proposition 1 of \cite{fasso1989composition} with $r=1$. 

\begin{proposition}
{\bf Exchange theorem: } Let $f$ be a real analytic function $f: U\times\torus^n 
\rightarrow\real$ extended to the domain $D_{\rho,\sigma}(U)$.  The equality
\begin{equation}\label{exchange}
f(\underline{A}',\underline{\varphi}') = 
\left(\exp(\mathcal{L}_\chi)f(\underline{A},\underline{\varphi})\right)
_{\underline{A}=\underline{A}',\underline{\varphi}=\underline{\varphi}'}
\end{equation}
holds, where $(\underline{A}',\underline{\varphi}')$ are given by the transformation 
(\ref{lietra}) and $(\underline{A}',\underline{\varphi}')\in 
D_{\rho-\delta,\sigma-\xi}(U)$.
\end{proposition}

See \cite{giorgilli2002notes} for the proof. In simple words, the exchange theorem 
implies that the result of a Lie series canonical transformation onto a function 
depending on $(\underline{A},\underline{\varphi})$ can be found by implementing 
the sequence of Poisson brackets of the exponential operator $\exp(\mathcal{L}_\chi)$ directly on the function $f$, and substituting, after this operation, the arguments 
$(\underline{A},\underline{\varphi})$ with $(\underline{A}',\underline{\varphi}')$.\\
\\
The above definitions allow us to establish an algorithm for the calculation of 
the sequence of canonical transformations (\ref{rseq}) using Lie series. The 
algorithm is obtained recursively by defining the $r$-th step as follows. Assume 
the Hamiltonian after $r-1$ normalization steps, denoted by $\H^{(r-1)}$, is in normal 
form up to the book-keeping order $r-1$:
\begin{equation}\label{hamrm1}
\H^{(r-1)}=Z_0+\epsilon Z_1+\ldots+\epsilon^{r-1}Z_{r-1} + 
\epsilon^{r} R^{(r-1)}_{r}+\epsilon^{r+1} R^{(r-1)}_{r+1}
+\epsilon^{r+2} R^{(r-1)}_{r+2}+\ldots
\end{equation}
Then, the $r$-th step Lie generating function $\chi_r$ and Hamiltonian $\H^{(r)}$ 
are computed as follows:\\
\\
\noindent
(i) split $R^{(r-1)}_{r}$ as $R^{(r-1)}_r = Z^{(r-1)}_r+h^{(r-1)}_r$, where 
$Z^{(r-1)}_r$ denotes the part of $R^{(r-1)}_{r}$ being in normal form; \\
\\
\noindent
(ii) compute $\chi_r$ as the solution of the homological equation
\begin{equation}\label{homo}
\{\underline{\omega}\cdot\underline{A},\chi_r\} + \epsilon^r h^{(r-1)}_r=0 ;
\end{equation}
\\
\noindent
(iii) compute the $r$-th step normalized Hamiltonian as $\H^{(r)} = 
\exp(\mathcal{L}_{\chi_r})\H^{(r-1)}$. This yields the Hamiltonian 
\begin{equation}\label{hamrm2}
\H^{(r)}=Z_0+\epsilon Z_1+\ldots+\epsilon^{r-1}Z_{r-1} + \epsilon^{r}Z_{r}+
\epsilon^{r+1} R^{(r)}_{r+1}+\epsilon^{r+2} R^{(r)}_{r+2}+\ldots
\end{equation}
where $Z_r = Z^{(r-1)}_r$. 

\begin{remark}
In the above algorithm, the notation $f^{(r)}$ implies a function depending on 
the canonical variables $(\underline{A}^{(r)},\underline{\varphi}^{(r)})$, 
which are connected to the original variables $(\underline{A},\underline{\varphi})$ 
via the composition of Lie series transformations
\begin{equation}\label{liecompo}
(\underline{A},\underline{\varphi}) = 
\exp(\mathcal{L}_{\chi_r})\circ\exp(\mathcal{L}_{\chi_{r-1}})\circ\ldots\circ
\exp(\mathcal{L}_{\chi_1})(\underline{A}^{(r)},\underline{\varphi}^{(r)}).
\end{equation}
For simplicity of notation, unless explicitly required in the sequel we do not 
write the superscripts in the canonical variables defined in every step, but only 
in the functions in which these variables are arguments of. 
\end{remark}

\begin{remark}
In the computer-algebraic implementation of the normalization algorithm, all functions 
are truncated up to a maximum book-keeping order, specified by computational 
restrictions. 
\end{remark}

\begin{remark}
The solution of the homological equation (\ref{homo}) is trivial when the functions 
$h^{(r-1)}_r$ are written in the Fourier representation
$$
h^{(r-1)}_r = \sum_{\underline{k}\in\mathbb{Z}^n} \tilde{h}^{(r-1)}_{r,\underline{k}}(\underline{A}) 
\exp(i\underline{k}\cdot\underline{\varphi}),
$$
which gives
$$
\chi_r = \epsilon^r \sum_{\underline{k}\in\mathbb{Z}^n} 
{i\tilde{h}^{(r-1)}_{r,\underline{k}}(\underline{A})\over
\underline{k}\cdot\underline{\omega}}\exp(i\underline{k}\cdot\underline{\varphi}).
$$
\end{remark}

\section{Stability of the semimajor axis in the $J_2$ model}\label{sec:stabilityJ2}

We will now implement the Hamiltonian normalization discussed in Section \ref{sec:NF} 
 to eliminate the short period terms (depending on the mean longitude $\lambda$) 
in the Hamiltonian (\ref{hj2exp}), leading to estimates on the long-term stability 
of the orbits' semimajor axis.

\subsection{Normal form}\label{ssec:nfj2} 
We express the Hamiltonian function in the form (\ref{bkeep}), choosing the 
book-keeping power equal to $s-2$, where $s$ is the index in the Hamiltonian expansion (\ref{hj2exp}), 
that is, collecting together at book-keeping order $s$ all polynomials 
$\mathcal{Z}_{s-2}$ and $\mathcal{P}_{s-2,k_1,k_2,k_3}$. Then
\begin{equation}
\begin{split}
&\mathcal{H}_{J_2}^{(0)}(\delta L, P, Q, \lambda,p,q)=\\&\quad
=\mathcal{H}_0(\delta L,P,Q)+\epsilon \mathcal{H}_1(\delta L, P, Q, 
\lambda,p,q)+\dots+\epsilon^N \mathcal{H}_N(\delta L, P, Q, \lambda,p,q)\ 
\end{split}
\end{equation}
with $\H_0 = n_*\delta L+\omega_{1*}P+\omega_{2*}Q$. The truncation order (in eccentricity and inclination) is  $N=15$. 

With reference to the algorithm of Subsection \ref{ssec: j2 gen}, normal form terms 
are specified as those non-depending on the mean longitude $\lambda$. After $M$ 
normalization steps, the Hamiltonian takes the form 
\begin{equation}\label{hamnew}
\H_{J_2}\equiv\mathcal{H}_{J_2}^{(M)}(\delta L, P, Q,
\lambda,p,q)=\mathcal{H}_{J_2,sec}^{(M)}(\delta L,P,Q,p,q)
+\mathcal{R}_{J_2}^{(M)}(\delta L,P,Q,\lambda, p,q)\ ,
\end{equation}
where (setting the book-keeping $\epsilon=1$)
\begin{equation*}
\begin{split}
\mathcal{H}_{J_2,sec}^{(M)} = Z_0+\ldots+Z_M,\\
\mathcal{R}_{J_2}^{(M)} = R^{(M)}_{M+1}+\ldots+R^{(M)}_{N}
\end{split}
\end{equation*}
with $Z_0=\H_0$. 

The term $\H_{J_2,sec}^{(M)}$ will be referred to as the `secular Hamiltonian' (not 
depending on the fast angle $\lambda$). On the other hand, the remainder $\R_{J_2}^{(M)}$ 
quantifies the difference between the true evolution of all canonical variables and 
the one induced by $\H_{J_2,sec}^{(M)}$. Since in (\ref{hamnew}) we can only compute a 
truncated remainder, we probe numerically that the finite sum of the leading terms 
in the remainder (up to order $N$) yields a remainder norm close to the limiting 
one (which corresponds to the limit $N\rightarrow\infty$). To this end, we take as 
maximum normalization order $M=N-3$, ensuring that at least the three first leading 
terms are included in the remainder (see \cite{steichen1997long}). Also, in estimating 
the size of the remainder through a suitable definition of the norm, we compute the 
\textit{sup} norm on a closed and bounded domain $\mathcal{D}\subset\real^2$:
\begin{equation}\label{supnorm}
\| f\|_{\infty, \mathcal{D}}=
\sup_{\substack{(e,i)\in\mathcal{D}\\(\lambda,p,q)\in\mathbb{T}^3}}
|f(e,i,\lambda,p,q)|\ .
\end{equation}
In practical computations, we can replace in \equ{supnorm} the supremum with the 
maximum, since the functions that we consider contain only positive powers of 
$\sqrt{1-e^2}$ and $\cos{i}$, hence, they are continuous on the domain $\mathcal{D}$. 

\subsection{Numerical results: stability of the semimajor axis}
\label{sec: stabest j2} 

Having fixed the procedure for the normal form and remainder computations, 
we proceed in deriving stability estimates based on the time variations of the value 
of the semimajor axis in the $J_2$ problem. Fixing a reference value $a_*$ of the 
semimajor axis, we assume that, at the time $t=0$, we have $L=L_*=\sqrt{\mu a_*}$, 
i.e. $\delta L=0$. Our aim is to estimate the fluctuations of $L$ as
functions of the orbital parameters $e$ and $i$.

The first question to settle is that, for every value of the
reference parameter $a_*$ we have to specify the range of values of
the variables $(e,i)$ for which the remainder $\mathcal{R}_{J_2}^{(M)}$ is small 
enough to represent
only a perturbation with respect to the dynamics determined by the
secular part. In applications, we compute the value of 
$\|\mathcal{R}_{J_2}^{(M)}\|_{\infty,\mathcal{D}}$ in the domain $(e,i)\in
\mathcal{D}=[0,0.15]\times[0,\pi/2]$, so that the inclination can take all
possible values; the eccentricity is instead taken in a reasonable interval,
where we can find almost all main Earth's satellites.

With reference to the Hamiltonian (\ref{hamnew}), if we consider the 
dynamics induced only by the
secular part, we obtain that
\[\frac{d}{dt}\delta L=-\frac{\partial \mathcal{H}_{J_2,sec}^{(M)}}{\partial \lambda}=0\ ,\]
which implies that $\delta L$ (hence $L$) is a constant of motion.
We remind that $\delta L$ is not the original Delaunay variable, but rather
the one obtained after $M$ normalization steps. If we
denote by $\delta L^{(0)}$ the original variable, then we have
\begin{equation}\label{deltaLfunctdeltaL0} \delta
L=\exp(-\mathcal{L}_{\chi^{(1)}}(\dots(\exp(-\mathcal{L}_{\chi^{(M)}}(\delta
L^{(0)})))))\ .
\end{equation}
To obtain $\delta L^{(0)}$ as a function of the new variable $\delta L$, we need 
to invert the transformation (\ref{deltaLfunctdeltaL0}); we observe that
\[\left(\exp(\mathcal{L}_{\chi})\right)^{-1}=\exp(-\mathcal{L}_{\chi})\ ,\]
implying
\[\delta L^{(0)}=\exp(\mathcal{L}_{\chi^{(M)}}
(\dots(\exp(\mathcal{L}_{\chi^{(1)}}(\delta L)))))\ .\]
Since we are dealing with \textit{near-identity} canonical
transformations, we realize that $\delta L^{(0)}$ is the sum of $\delta L$
and short period (small) variations {which do not affect its stability}.

If we consider the full Hamiltonian in (\ref{hamnew}), then $L$ is not 
constant anymore because of the dependence of $\R_{J_2}^{(M)}$ on $\lambda$. Using again 
Hamilton's equations, we see that
\[\frac{d}{dt}L=\frac{d}{dt}(\delta L+L_*)=\frac{d}{dt}\delta L
=-\frac{\partial \mathcal{H}_{J_2}}{\partial \lambda}
=-\frac{\partial \R_{J_2}^{(M)}}{\partial \lambda}\ .\]
Then, for every set of values, say 
$(e^*,i^*,\lambda^*,p^*,q^*)\in\mathcal{D}\times\mathbb{T}^3$, 
we obtain
\[\bigg|\frac{d}{dt}L(e^*,i^*,\lambda^*,p^*,q^*)
\bigg|\leq\sup_{\substack{(e,i)\in\mathcal{D}\\(\lambda,p,q)\in\mathbb{T}^3 }}\bigg|
\frac{d}{dt}L(e,i,\lambda,p,q)\bigg|=\bigg|\bigg|\frac{\partial
\R_{J_2}^{(M)}}{\partial\lambda}\bigg|\bigg|_{\infty, \mathcal{D}}\ .\]
Let $L(e,i,\lambda,p,q;T)$ be the value at time $t=T$. To estimate its distance 
from the equilibrium point $L_*$, we can use the mean value theorem which gives
\begin{equation}\label{stabest}
\vert L(e,i,\lambda,p,q;T)-L_*\vert\leq\|L(e,i,\lambda,p,q;T)-L_*\|_{\infty,\mathcal{D}}
\leq\bigg|\bigg|\frac{dL}{dt}\bigg|\bigg|_{\infty,\mathcal{D}}T\ .
\end{equation}
Requiring that the right hand side of (\ref{stabest}) is of order of unity, 
then the stability time $T$ becomes order of 
$O\left(1/\| dL/dt\|_{\infty,\mathcal{D}}\right)$. Let us fix a constant 
value $\Delta L$ and suppose that we want to estimate the minimal time $T_1$ 
up to which the variation of $L(e,i,\lambda,p,q;T)$ stays bounded by $\Delta L$:
\[\| L(e,i,\lambda,p,q;T)-L_*\|_{\infty,\mathcal{D}}\leq\Delta L\ .\]
Using (\ref{stabest}) we obtain that $T_1$ is given by
\begin{equation}\label{T1}
T_1\geq\frac{\Delta L}{\| dL/dt\|_{\infty,\mathcal{D}}}\ .
\end{equation}
Equation (\ref{T1}) can be used to derive the stability time of the semimajor axis $a$: recalling that, in general, $L=\sqrt{\mu a}$, one has that $\Delta L=\Delta a/2\sqrt{\mu/ a}$, which allows to obtain a lower bound for the stability time of $a$ given by 
\begin{equation}\label{T2}
T_2=\frac{1}{2}\sqrt{\frac{\mu}{a_*}}\frac{\Delta a}{\| dL/dt\|_{\infty,\mathcal{D}}}.
\end{equation}
This estimate will be used in Section~\ref{ssec: stabest ls num} to obtain 
results on the stability time at different altitudese; in particular, $\Delta a$ is set to be equal to $0.1$ $R_E$. 

To check that the norm $\|\R_{J_2}^{(M)}\|_{\infty,\mathcal{D}}$ is \textit{small} in
the domain $\mathcal{D}=[0,0.15]\times[0,\pi/2]$, we compute its value by 
taking a set of samples for the reference value of the semimajor axis $a_*$, 
that correspond to different distances from the Earth's center (the radius 
of the Earth is $R_E=6378.14$~km). Precisely, we consider the following 
semimajor axes:
\begin{itemize}
\item $a_*^{(1)}=(42164$ $km)/R_E$: the reference value for GEO satellites;
\item $a_*^{(2)}=(26560$ $ km)/R_E$: the reference value for GPS satellites;
\item $a_*^{(3)}=(8524.75 $ $ km)/R_E$: an intermediate value in terms of the 
altitude;
\item $a_*^{(4)}=(7258.69$ $ km)/R_E$: very close to the Earth's surface.
We remark that in this case the results obtained are not very relevant 
from a practical point of view, because the effect of the atmosphere 
becomes important.
\end{itemize}
\begin{table}
    \caption{Estimates of $\| \R_{J_2}^{(M)}\|_{\infty,\mathcal{D}}$ for different 
    values of $a_*$ in the $J_2$ model.}
    \label{tabellaresto}
    \begin{center}
        \begin{tabular}{|c|c|c|}
            \hline
            Semimajor axis & $a_*$ & $\| \R_{J_2}^{(M)}\|_{\infty,\mathcal{D}}$ \\
            \hline
            $42164\quad km$ & $6.6107$ & $1.28967\cdot 10^{-11}$  \\
            \hline
            $ 26560\quad km$ & $4.16422$ & $1.60737\cdot10^{-10}$  \\
            \hline
            $8524.75\quad km$ & $1.33656$ & $6.26588\cdot10^{-8}$  \\
            \hline
            $7258.69\quad km$ & $1.13806$ & $1.43864\cdot10^{-7}$  \\
            \hline
    \end{tabular}   \end{center}
\end{table}
Table \ref{tabellaresto} shows the values of $\|\R_{J_2}^{(M)}\|_{\infty,\mathcal{D}}$ 
computed for the above values of $a_*$ and for $J_2=1.084\cdot10^{-3}$, 
namely the real value of the coefficient for the Earth. As we can see,
$\|\R_{J_2}^{(M)}\|_{\infty,\mathcal{D}}$ is typically very small for all values of 
$a_*$: this confirms that for the $J_2$ problem it is reasonable to take
the domain in eccentricity and inclination as 
$\mathcal{D}=[0,0.15]\times[0,\pi/2]$.

\begin{table}
    \caption{Estimates of $\| dL/dt\|_{\infty,\mathcal{D}}$ for different 
    values of $a_*$ in the $J_2$ model.}
    \label{tab:tabelladeL}
    \begin{center}
        \begin{tabular}{|c|c|c|}
            \hline
            Semimajor axis & $a_*$ & $\| dL/dt\|_{\infty,\mathcal{D}}$ \\
            \hline
            $42164\quad km$ & $6.6107$ & $2.7216\cdot10^{-10}$  \\
            \hline
            $ 26560\quad km$ & $4.16422$ & $6.66832\cdot10^{-10}$  \\
            \hline
            $8524.75\quad km$ & $1.33656$ & $1.63251\cdot10^{-7}$  \\
            \hline
            $7258.69\quad km$ & $1.13806$ & $3.4383\cdot10^{-7}$  \\
            \hline
    \end{tabular}   \end{center}
\end{table}
Table \ref{tab:tabelladeL} provides the results 
for the estimate of $\|dL/dt\|_{\infty,\mathcal{D}}$, which show that, using 
\equ{T1} with $\Delta L$ equal for all the considered distances $a_*$, the stability time for $L$ increases with the altitude.

\begin{figure}[h]
\vglue6cm
\hglue-8cm
    \includegraphics[width=0.1truecm,height=0.1truecm]{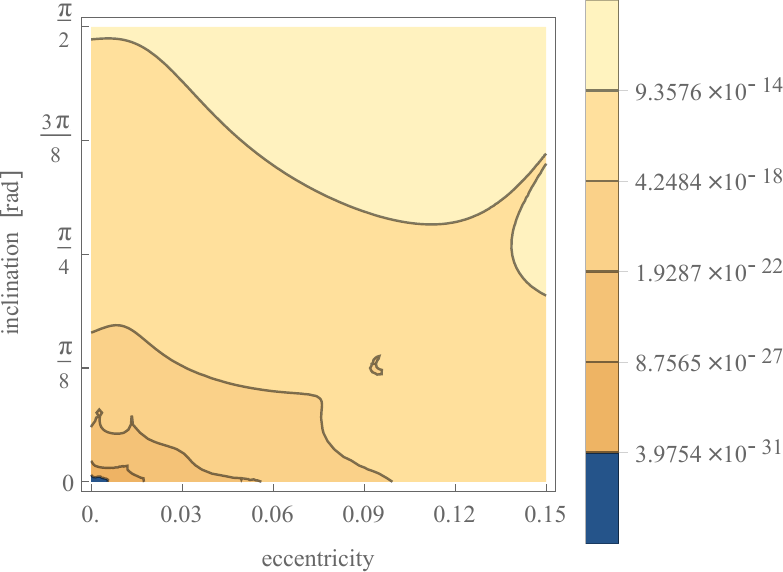}
\hglue8.5cm
    \includegraphics[width=0.1truecm,height=0.1truecm]{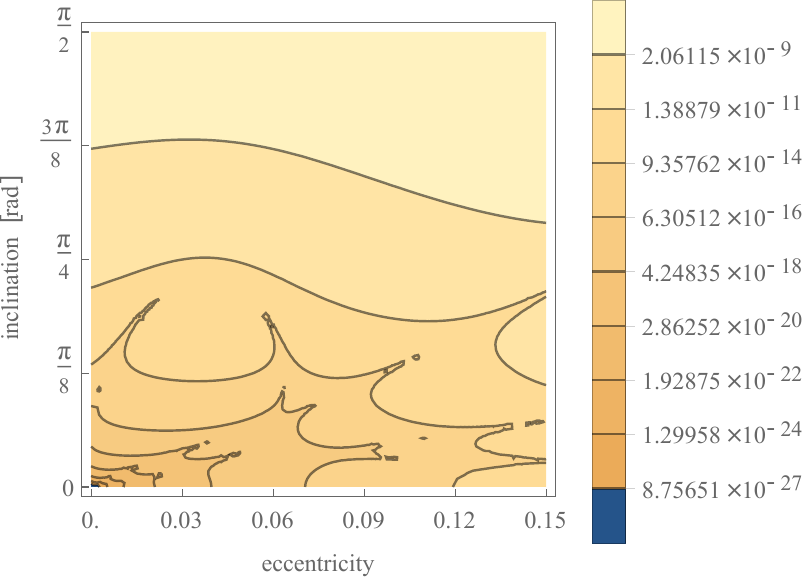}
\hglue-1cm
\caption{Plots of $\|\R_{J_2}^{(M)}\|_{\infty,\mathcal{D}}$ for
$(e,i)\in\mathcal{D}$: $a_*=a_*^{(1)}$ (left) and
$a_*=a_*^{(4)}$ (right) in the $J_2$ model.}
    \label{fig:normerestogeoene}
\end{figure}

\begin{figure}[h]
\centering
\vglue7cm
\hglue-8.5cm
    \includegraphics[width=0.1cm,height=0.1cm]{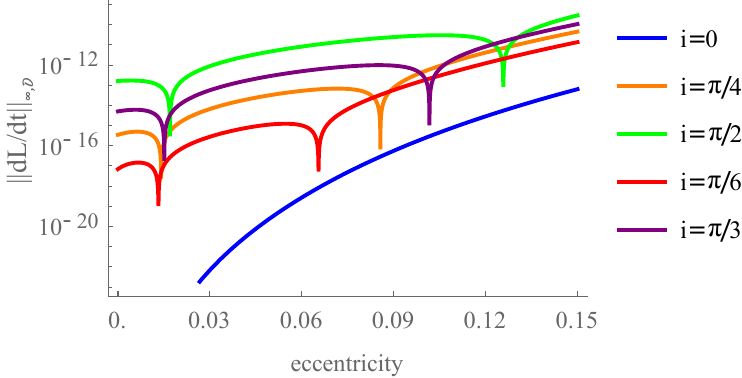}
\hglue8.5cm
    \includegraphics[width=0.1cm,height=0.1cm]{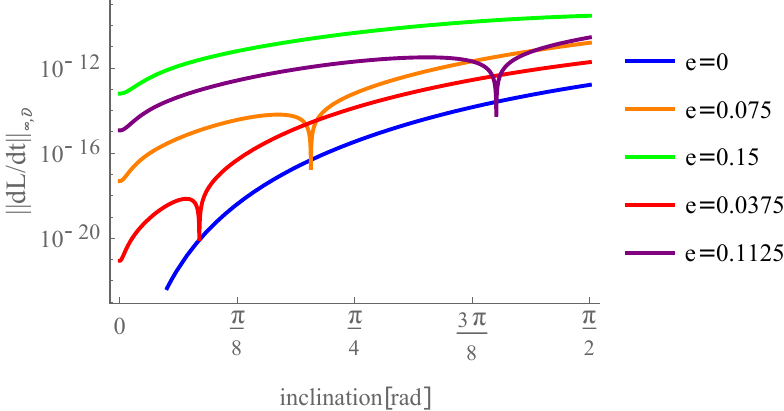}
    \caption{Stability estimates for the GEO case (with $a=a_*^{(1)}$) in the $J_2$ model. Left: 
    plot of $||dL/dt||_{\infty,\mathcal{D}}$ as a function of $e$
    for fixed values of $i$. Right: plot of $||dL/dt||_{\infty,\mathcal{D}}$ 
    as a function of $i$ for fixed values of $e$.
    }
    \label{fig:normastabgeo}
\end{figure}

Figure~\ref{fig:normerestogeoene} shows the logarithmic plot of 
$\|\R_{J_2}^{(M)}\|_{\infty,\mathcal{D}}$ in the limit cases $a_*=a_*^{(1)}$ and $a_*=a_*^{(4)}$. 
The plots show that the remainder decreases as one gets farther from the Earth 
and it becomes larger when increasing the eccentricity and inclination.

\begin{figure}[h]
\centering
\vglue4cm
\hglue-7.5cm
    \includegraphics[width=1cm,height=1cm]{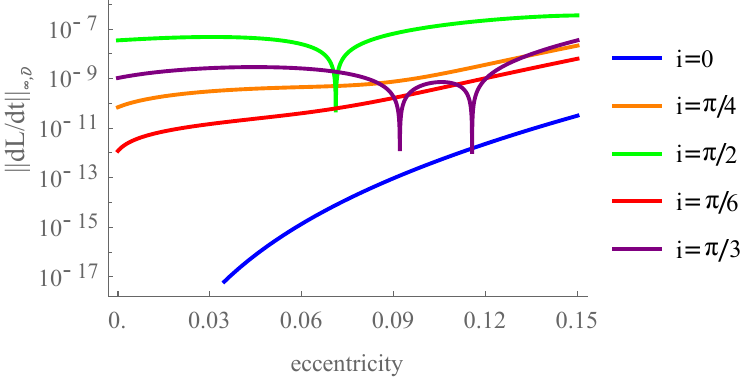}
\hglue7.5cm
    \includegraphics[width=1cm,height=1cm]{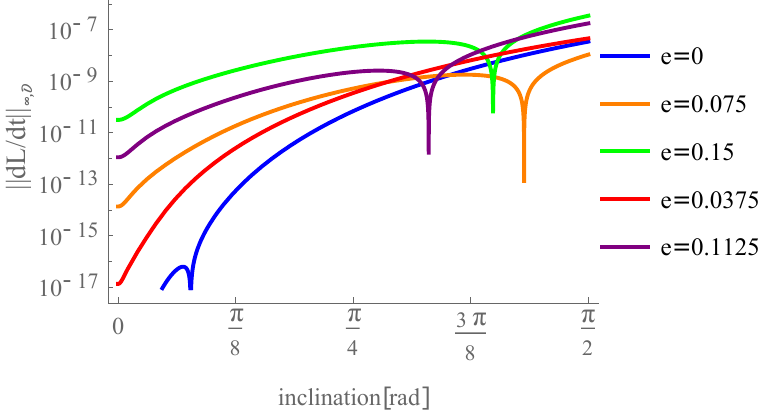}
    \caption{Stability estimates for the near-Earth case (with $a=a_*^{(4)}$) in the $J_2$ model. 
    Left: plot of $||dL/dt||_{\infty,\mathcal{D}}$ as a function of $e$
    for fixed values of $i$. Right: plot of $||dL/dt||_{\infty,\mathcal{D}}$ 
    as a function of $i$ for fixed values of $e$.
    }
    \label{fig:normastabne}
\end{figure}

Figures~\ref{fig:normastabgeo} and \ref{fig:normastabne} refer, respectively, 
to $a=a_*^{(1)}$ and $a=a_*^{(4)}$; the left plots provide the graph of 
$||dL/dt||_{\infty,\mathcal{D}}$ as a function of the eccentricity for fixed
values of the inclination, while the right plots give the norm as a function 
of the inclination for fixed values of the eccentricity. We notice that the norms 
tend to decrease when the eccentricity and the inclination are smaller,
although the effect is more evident in the GEO region than closer to the Earth.

We now examine how the stability time changes as a function of the semimajor 
axis $a_*$: in this case, we consider 1000 values for $a_*$ uniformly distributed 
from $a_{in}=1.15679$ (corresponding to an altitude of $1000$ $ km$, which we take as the first reference value, although in this region weak dissipative effects are possibily affecting the dynamics) to 
$a_f=16.6786$ (corresponding to an altitude of $10^5$ $ km$), using Eq. (\ref{T2}) with $\Delta a=0.1$ $R_E$. 
Figure \ref{fig:stabilityinfunctionofa} confirms that the stability time 
increases with the altitude also in the case of the semiajor axis. Indeed, while for $a_*=a_{in}$ we can ensure the 
stability of the semimajor axis for a period of the order of years, in the 
case $a_*=a_{f}$ we have a stability time of the order of $10^4$ years. From 
an analytical point of view, this behaviour of the stability time can be 
explained by the fact that, for higher distances, our model can be approximated 
by Kepler's problem in which the semimajor axis is constant.

\begin{figure}[h]
\centering
\vglue6cm
\hglue-10cm
    \includegraphics[width=1cm,height=1cm]{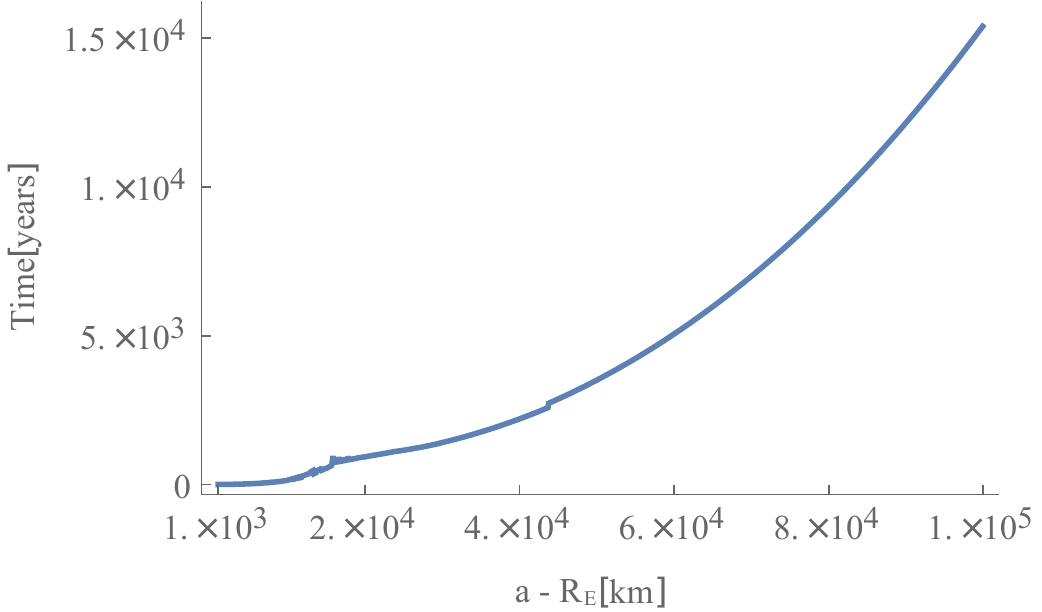}
    \caption{Stability time in the $J_2$ model for $a\in[a_{in},a_{f}]$ (see the text for the 
    definition of $a_{in}$, $a_{f}$) allowing a variation of $0.1$ $R_E$. }
    \label{fig:stabilityinfunctionofa}
\end{figure}

\section{Secular stability in the geolunisolar model}\label{sec: stabest ls}
Using the $J_2$ model, we have demostrated how the stability of the semimajor 
axis can be established against {\it short-period} perturbations (depending 
on the satellite's mean anomaly). In this section, we focus, instead, on 
the long-term variations in the eccentricity and inclination of the satellite's orbit, for orbits close to circular ($e<0.1$ $rad$) and with small inclination 
$(|i|<0.1)$. One can easily verify that, within the geolunisolar problem 
(Hamiltonian $\H_{gls, sec}^{\leq N}$, see Eq. (\ref{hamavfinaaN})), the phase-space manifold $e=0$, corresponding 
to $I_2=0$, constitutes an invariant manifold of the flow, implying that 
circular orbits remain so for infinitely long times independently of their 
variations in inclination and longitude of the node. On the other hand, 
for $e$ small, but non-zero, long-term variations of both the 
eccentricity and inclination can occur on timescales given by the inverse 
of the frequencies $\nu_1$ and $\nu_2$ (Eq. (\ref{hamavfinpol})). Since 
$\nu_1\simeq \nu_2 \simeq \frac{3}{2}\frac{\sqrt{\mathcal{G}M_E}R_E^2J_2}{a^{7/2}}$, 
the secular timescale is of order of
$T_{sec}=\mathcal{O}\left((a/R_E)^2J_2^{-1}\right)T_{short}$, where $T_{short}$ is the characteristic time of the frequency associated to the fast angle. Since $J_2\approx 
10^{-3}$, the short and long periods are separated by three orders of magnitudes, 
a fact which justifies altogether the simple averaging over mean anomalies which leads to the model of departure $\H_{gls, sec}^{\leq N}$ for the analysis of the secular 
stability. On the other hand, the fact that $\nu_1\simeq\nu_2$ 
implies that, near the equator (or, more precisely, for orbits near the Laplace 
plane, see Section \ref{ssec: forced}), the eccentricity and inclination have coupled 
variations (the so-called 'Kozai-Lidov' mechanism). This fact implies that, 
close to the Laplace plane, the term `secular stability' cannot mean the 
long-term preservation of the eccentricity and inclination one independently 
of the other, but only the approximate preservation of the combination 
$\mathcal{I}\approx e^2+i^2$ (see below for exact expressions) known as the 
Kozai-Lidov integral. The normal form construction and remainder estimates 
in the present section reflect these basic properties of the dynamics. 

\subsection{Normal form}\label{ssec: NFjls}
Starting with the model $\H_{gls, sec}^{\leq N}$ given in Eq. (\ref{hamavfinaaN}), the 
construction of the normal form proceeds with the algorithm described in 
Section \ref{sec:NF} and the following settings: \\
\\
i) The book-keeping rule 
(exponent $s$ in Eq. (\ref{bkeep})) is set as $s=s_1+s_2-2$, where $s_1$ 
and $s_2$ are the exponents appearing in Eq. (\ref{hamavfinaaN}). \\
ii) The resonant module (Eq. (\ref{nfres}), case 3 of subsection 3.1) is 
set as: 
$$
\mathcal{M}:=\{(k_1,k_2)\in\mathbb{Z}^2: k_1+k_2=0\}
$$
where $k_1,k_2$ are the integers specifying each Fourier harmonic in 
Eq. (\ref{hamavfinaaN}).\\
iii) The maximum truncation order is set to $N=15$, while the maximum 
normalization order is set to $M=12$. 

With the following settings, the Hamiltonian after $r$ normalization steps, 
where $r$ can take the values $r=1,2,...M$, resumes the form:
\begin{equation}\label{Hglsr}
\mathcal{H}_{gls, sec}^{(r)}(I_1,I_2, \phi_1, \phi_2)=\mathcal{Z}_{gls, sec}^{(r)}(I_1, I_2)
+\mathcal{Z}_{gls, res}^{(r)}(I_1, I_2, \phi_1-\phi_2)+\R_{gls}^{(r)}(I_1,I_2,\phi_1,\phi_2)~~.
\end{equation}
The term $\mathcal{Z}_{gls, sec}^{(r)}(I_1, I_2)$, hereafter called the \textit{secular part},
contains all terms independent of the angles (corresponding to the choice $k_1=k_2=0$ in the resonant module). The dynamics of this term implies separate 
preservation of the eccentricity and inclination (the latter around the Laplace 
plane). Instead, $\mathcal{Z}_{gls, res}^{(r)}(I_1,I_2, \phi_1-\phi_2)$, called the 
\textit{resonant part} of the normal form, collects all normal form terms 
depending on the resonant angle $\phi_1-\phi_2$. Finally, $\R_{gls}^{(r)}(I_1,I_2, \phi_1,\phi_2)$ is the \textit{remainder term}, which contains non-normalized terms of book-keeping 
orders $s=r+1,\ldots,N$. After $M$ normalization steps, we obtain the final geolunisolar Hamiltonian $\H_{gls}\equiv\H_{gls,sec}^{(M)}$.

We now look at the dynamics induced by the sum of secular and resonant parts:
\begin{equation*}
\mathcal{H}_{norm}(I_1,I_2,\phi_1,\phi_2)=\mathcal{Z}_{gls, sec}^{(M)}(I_1,I_2)+\mathcal{Z}_{gls, res}^{(M)}(I_1,I_2,\phi_1-\phi_2),
\end{equation*}
called, altogether, the \textit{resonant normal form} $\mathcal{H}_{norm}$ (for 
simplicity, we drop the dependence on the normalization order $r$ from the 
notation). The quantity $I_1+I_2$ is a first integral for the dynamics induced by $\mathcal{H}_{norm}$,
which implies that the vertical component of the angular momentum, which coincides
with $\Theta$, is preserved\footnote{For the
	$J_2$ model the preservation of the vertical component of the angular
	momentum is a direct consequence of the axisimmetry of the
	truncated geopotential. For the geolunisolar model, the addition of
	the external attractions breaks this symmetry. However, the
	preservation of this quantity turns to be still true for the
	Hamiltonian $\mathcal{H}_{norm}$.}. Given that $L$ is constant,
say $L=L_*=\sqrt{\mu_E a_*}$, the quantity
\begin{equation}\label{I1+I2}
I_1+I_2=L_*-L_*\sqrt{1-e^2}(1-\cos{i})
\end{equation}
is a first integral and, as a consequence, the quantity
\begin{equation}\label{glsint}
\mathcal{I}(e,i)=1-\sqrt{1-e^2}(1-\cos{i})
\end{equation}
is constant for the dynamics induced by the normal form. This means that $e$ and $i$ can change only in such a way that the value of $\mathcal{I}(e,i)$ remains constant.

The fact that the presence of resonant first integrals determines a \textit{locking} in the values of $e$ and $i$ is at
the basis of the so-called \textit{Lidov-Kozai effect} (\cite{Lidov,Kozai}), which is common, in a wide range of resonant combinations,
in many models of Celestial Mechanics.

\subsection{Remainder and stability estimates}
As already mentioned in Section \ref{sec: stabest j2},  we need to guarantee that the remainder is \textit{small}
with respect to the normal part; we denote again by $\mathcal{D}$ the domain over which the norm
$\|\R_{gls}^{(M)}\|_{\infty,\mathcal{D}}$ is computed where, for a function $f=f(e,i,\phi_1,\phi_2)$, the norm of $f$ is defined as
$$
\| f\|_{\infty,\mathcal{D}}=\sup_{(e,i)\in\mathcal{D},(\phi_1,\phi_2)\in\torus^2} |f(e,i,\phi_1,\phi_2)|\ .
$$
There exists an optimal value of $M$ that minimizes the estimate of the remainder's norm,
as shown in Section~\ref{ssec: stabest ls num} for GEO orbits.

Since $I_1+I_2$ is a first integral for $\mathcal{H}_{norm}$, we have
that
\begin{equation*}
\{I_1+I_2, \mathcal{H}_{norm}\}=0.
\end{equation*}
To evaluate the stability of $\mathcal{I}(e,i)$, we use the relation:
\begin{equation*}
\frac{d}{dt}(I_1+I_2)=\{I_1+I_2, \mathcal{H}_{gls}\}=\{I_1+I_2, \R_{gls}^{(M)}\};
\end{equation*}
then, for every $(e^*,i^*,\phi_1^*,\phi_2^*)\in \mathcal{D}\times\mathbb{T}^2$, we have the following estimate:
\begin{equation*}
\bigg\vert
\frac{d}{dt}(I_1+I_2)(e^*,i^*,\phi_1^*,\phi_2^*)\bigg\vert\leq\sup_{\substack{(e,i)\in\mathcal{D}\\(\phi_1,\phi_2)\in\mathbb{T}^2
}}\bigg\vert
\frac{d}{dt}(I_1+I_2)(e,i,\phi_1,\phi_2)\bigg\vert\leq\|\{I_1+I_2,\R_{gls}^{(M)}\}\|_{\infty,
	\mathcal{D}}\ .
\end{equation*}
Let us now consider an orbit with initial point $(I_{1,0}, I_{2,0})$ such that the
corresponding eccentricity and inclination belong to $\mathcal{D}$; consider
its evolution up to $t=T$. Using the mean value theorem, we have that
\begin{equation}\label{lsstab}
\| (I_1(T)+I_2(T))-(I_{1,0}+I_{2,0})\|\leq\|\{I_1+I_2,\R_{gls}^{(M)}\}\|_{\infty, \mathcal{D}}\ T.
\end{equation}
Setting $\Gamma$ to be the maximum value for the variation of
$I_1+I_2$ in time, let us denote by $\widetilde{T}$ the minimum time such
that for every $T\leq\tilde{T}$
\begin{equation*}
\| (I_1(T)+I_2(T))-(I_{1,0}+I_{2,0})\|\leq\Gamma\ .
\end{equation*}
From (\ref{lsstab}), we have
\begin{equation*}
\widetilde{T}\geq\frac{\Gamma}{\|\{I_1+I_2, \R_{gls}^{(M)}\}\|_{\infty,
		\mathcal{D}}};
\end{equation*}
then we can use the value of $T$ as $T=\Gamma/\|\{I_1+I_2,\R_{gls}^{(M)}\}\|_{\infty, \mathcal{D}}$, which gives an estimate for the
stability time of $I_1+I_2$ and, consequently, of $\mathcal{I}(e,i)$. The stability results for the quantity $\mathcal{I}$ can be translated in terms of the orbital elements $(e,i)$ as follows: in view of (\ref{glsint}), for small values of $e$ and $i$ we find
\begin{equation}
\mathcal{I}\simeq L_*\frac{e^2+i^2}{2}, 
\end{equation}
hence, if we consider the variations of $\mathcal{I}$, $e$ and $i$, they are connected by the relation
\begin{equation}
\frac{\Delta\mathcal{I}}{\mathcal{I}}\simeq 2\frac{e\Delta e+i\Delta i}{e^2+i^2}. 
\end{equation}
For the limit case of $e$ or $i$ fixed and small, one find 
\begin{equation}
\frac{\Delta \mathcal{I}}{\mathcal{I}}\simeq 2\frac{\Delta e}{e}\simeq 2\frac{\Delta i}{i}, 
\end{equation}
then the relative variation of $\mathcal{I}$ (and, as a consequence, of $I_1+I_2$) is proportional to the relative variations of the orbital elements by a factor 2. \\
To make the stability results for the geolunisolar model consistent with the ones obtained in Section \ref{sec: stabest j2} for the $J_2$ model, in Section \ref{ssec: stabest ls num} we set 
\begin{equation}
\Gamma=\frac{1}{2}\sqrt{\frac{\mu}{a_*}}0.1,
\end{equation}
namely, recalling that $\Delta L=\Delta a/2\sqrt{\mu/a}$ and $\Delta a=0.1$, the maximal variation of $I_1+I_2$ in the geolunisolar model is equal to the maximal variation allowed  for the action $L$ in the $J_2$ model. 

\subsection{Numerical results for the geolunisolar model}\label{ssec: stabest ls num}
For the geolunisolar model, we take the domain $(e,i)\in \mathcal{D}=[0,0.1]\times[0,0.1]$ around the forced eccentricity (which is always zero) and the forced inclination (which depends on the chosen altitude).

Since the stability results strongly depend on the distance from
the Earth, we select five different altitudes, that correspond to
cases of interest for the satellite's problem:
\begin{itemize}
    \item $h^{(1)}=3000$ $km$, above the atmosphere;
    \item $h^{(2)}=20000$ $km$, that is in MEO region;
    \item $h^{(3)}=35786$ $km$, the altitude of GEO orbits;
    \item $h^{(4)}=50000$ $km$, corresponding to far objects;
    \item $h^{(5)}=100000$ $km$, that corresponds to objects which are very far from the Earth's surface.
\end{itemize}
The value of the remainder's norm depends on the altitude of the
orbit: in particular, we can state that the stability time
decreases as the altitude increases.

\begin{table}
    \caption{Estimate of $\|\R_{gls}^{(M)}\|_{\infty, \mathcal{D}}$, with $M=12$, in the geolunisolar model with $\mathcal{D}=[0,0.1]\times[0,0.1]$
    for different altitudes.}
    \label{conv in 0.1 0.1}
    \begin{center}
        \begin{tabular}{|c|c|}
            \hline
            \textbf{Altitude} &$\| R^{(M)}\|_{\infty, \mathcal{D}}$\\
            \hline
            $3000$ $km$ & $3.74442\cdot10^{-16}$ \\
            \hline
            $20000$ $km$ & $1.82777\cdot10^{-15}$ \\
            \hline
            $35790$ $km$ & $8.7787\cdot10^{-15}$ \\
            \hline
            $50000$ $km$ & $4.97867\cdot10^{-12}$ \\
            \hline
            $100000$ $km$ &$ 1.64614\cdot10^{-9} $ \\
            \hline
    \end{tabular}   \end{center}
\end{table}

\begin{figure}[!t]
\centering
\vglue5cm
\hglue-6cm
    \includegraphics[width=1cm,height=1cm]{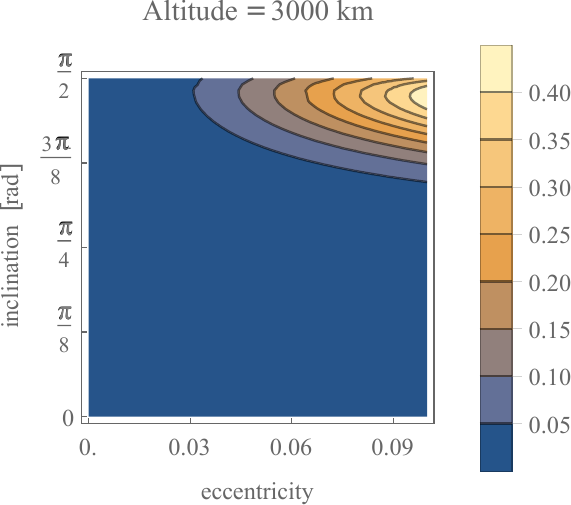}
\hglue7.5cm
    \includegraphics[width=1cm,height=1cm]{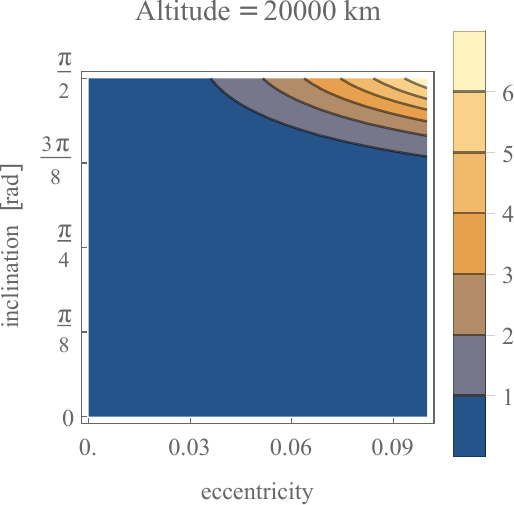}
\vglue4.7cm
\hglue-6cm
    \includegraphics[width=1cm,height=1cm]{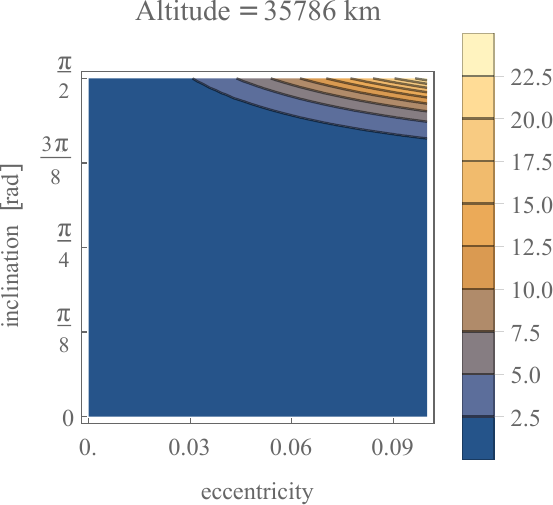}
\hglue7.5cm
    \includegraphics[width=1cm,height=1cm]{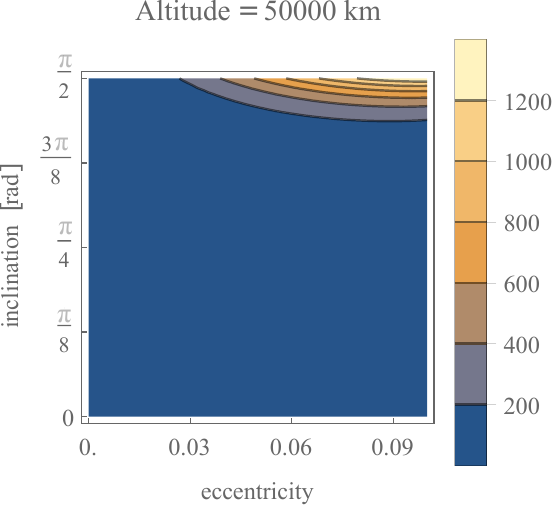}
\vglue4.7cm
\hglue-6cm
 \includegraphics[width=1cm,height=1cm]{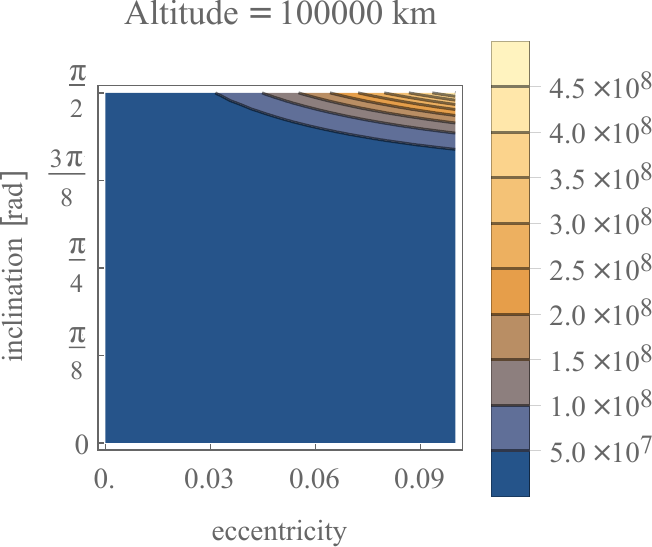}

\caption{Remainder's norm for the geolunisolar model in the domain $\mathcal{D}'=[0,0.1]\times[0,\pi/2]$ for $h^{(i)}$, $i=1,\dots, 5$
    (see the text for the definition of $h^{(i)}$).}
    \label{fig:normerestols}
\end{figure}

Table \ref{conv in 0.1 0.1} provides the value of $\|\R_{gls}^{(M)}\|_{\infty,\mathcal{D}}$
as a function of the altitude, showing a significant worsening for altitudes after the GEO region.

Figure \ref{fig:normerestols} shows the behaviour of the remainder's norm as a function of $(e,i)$ in
the bigger domain $\mathcal{D}'=[0,0.1]\times[0,\pi/2]$: as we can see, in almost all
cases the domain $\mathcal{D}'$ is too large to ensure the smallness of $\|\R_{gls}^{(M)}\|_{\infty,\mathcal{D}'}$.
Moreover, the magnitude of $\|\R_{gls}^{(M)}\|_{\infty,\mathcal{D}'}$ increases significantly with the altitude.
We can easily notice that the value of $\|\R_{gls}^{(M)}\|_{\infty,\mathcal{D}'}$ is strongly
dependent on the inclination: using this fact, we can detect a
value of $i$, denoted by $i_{crit}$, which is the minimum value of
the inclination for which $\|\R_{gls}^{(M)}\|_{\infty,\mathcal{D}'}$ is of the order of unity.
\begin{table}
    \caption{Value of $i_{crit}$ as a function of different values of the altitude for the geolunisolar model.}
    \label{icrit ls}
    \begin{center}
        \begin{tabular}{|c|c|c|}
            \hline
            \textbf{Altitude} &$i_{crit}(deg)$\\
            \hline
            $3000$ $km$ &$90^\circ$\\
            \hline
            $20000$ $km$ &$68.75^\circ$\\
            \hline
            $35790$ $km$ &$68.75^\circ$\\
            \hline
            $50000$ $km$ &$60.73^\circ$\\
            \hline
            $100000$ $km$ &$24.06^\circ$\\
            \hline
    \end{tabular}   \end{center}
\end{table}
Table \ref{icrit ls} shows the computed values of $i_{crit}$ (converted in degrees) for
the considered altitudes: we can notice that the
\textit{smallness} domain shrinks substantially between $50000$
$km$ and $100000$ $km$; in any case, we can see that for every value ofthe considered altitudes the domain $\mathcal{D}=[0,0.1]\times[0,0.1]$ is contained in the smallness domain of $\R_{gls}^{(M)}$. 

As mentioned in Section \ref{sec: stabest ls}, the remainder's
norm depends on the normalization order
$M$. Although the norm does not converge to zero if
$M$ tends to infinity, there is a value
of $M$, called the \textit{optimal normalization order}, say
$M_{opt}$, for which the norm of the remainder is minimal. Typically,
this optimal value is greater than the order of the Taylor
expansions of the numerically computed functions, and the
estimates for the remainder is so good that there is no reason to
push further the order of the expansion; for example, this is the case for
the normalized Hamiltonian function which describes the geolunisolar
problem computed for the GEO altitude.

\begin{figure}[!h]
\centering
\vglue4cm
\hglue-8cm
\includegraphics[width=1cm,height=1cm]{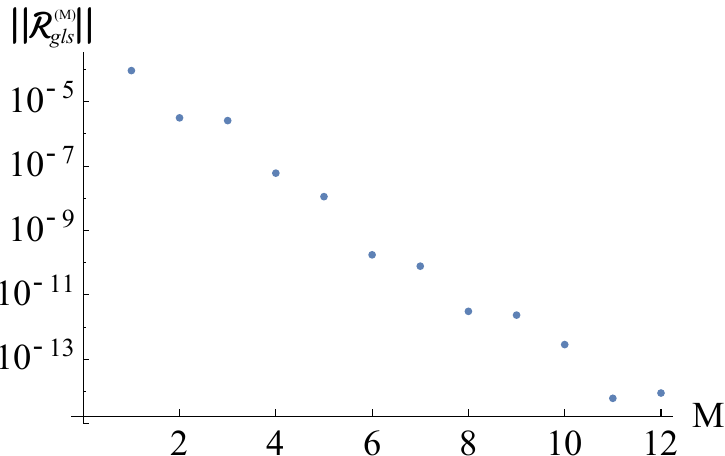}\\
\caption{Estimate of $\|\R_{gls}^{(M)}\|_{\infty,
\mathcal{D}}$ as a function of the normalization order $M$ for
$h=h^{(3)}$ (GEO distance) in the domain
$\mathcal{D}=[0,0.1]\times[0,0.1]$ for the geolunisolar model.}
    \label{fig:moptgeo}
\end{figure}

As we can see from Figure \ref{fig:moptgeo}, the optimal normalization order is greater than or equal to $11$. Since the values of $\|\R_{gls}^{(M)}\|_{\infty,\mathcal{D}}$ are very close for these values of $M$, for simplicity we made our estimates for the last computed normalization order $M=12$. 

Once obtained the smallness of $\|\R_{gls}^{(M)}\|_{\infty,\mathcal{D}}$ in
$\mathcal{D}$, we proceed to compute the stability time for the quantity
$\mathcal{I}(e,i)=1-\sqrt{1-e^2}(1-\cos{i})$.

\begin{table}
    \caption{Stability times in years for different altitudes in the domain $(e,i)\in \mathcal{D}=[0,0.1]\times[0,0.1]$ for the geolunisolar model.}
    \label{stab time}
    \begin{center}
        \begin{tabular}{|c|c|}
            \hline
            \textbf{Altitude} &\textbf{Stability time in} $\mathcal{D}$\\
            \hline
            $3000$ $km$ & $3.86102\cdot10^{14}$ \\
            \hline
            $20000$ $km$ & $4.49464\cdot10^{14}$\\
            \hline
            $35790$ $km$ & $1.17054\cdot10^{14}$\\
            \hline
            $50000$ $km$ & $ 2.11051\cdot10^{10}$\\
            \hline
            $100000$ $km$ &$1.21928\cdot10^7 $\\
            \hline
    \end{tabular}   \end{center}
\end{table}

As we can see from Table \ref{stab time},
the stability times are extremely long: this fact depends on
the model we considered, with the Lunar orbit in the ecliptic plane without precession effects.
However, we can notice a relevant decrease in the stability time for distances greater than GEO.
This behaviour is opposite to that of the $J_2$ model where the stability time was increasing
with the altitude (see Figure~\ref{fig:stabilityinfunctionofa}). In fact, at low altitudes the $J_2$ model is
strongly affected by the Keplerian part and the geopotential, while the geolunisolar
model takes into account both the inner effect due to the Earth and the outer effect due
to Moon and Sun.

As a final remark, to show the importance of taking the right domain in eccentricity
and inclination, let us assume $h=h^{(5)}$ and consider the
domain $(e,i)\in B=[0,0.1]\times[0,0.5]$, which is larger than the
convergence domain $[0,0.1]\times[0,i_{crit}]$ (see Table
\ref{icrit ls}). If we compute the stability time in the enlarged
domain $B$, we obtain just the value $T=0.00164$ years.

\section{Non-degeneracy conditions}\label{sec: nondeg}
Beside the orbital stability obtained using the Birkhoff normal form as in the previous Sections,
analytical estimates of the stability time can be obtained by the outstanding theorem developed by Nekhoroshev
(\cite{nekhoroshev1977exponential}). Under suitable assumptions, the theorem gives a confinement of the action variables for exponentially
long times. In particular, the Hamiltonian must satisfy a non-degeneracy condition which, in the original formulation,
is called \sl steepness \rm condition. The definition of the steepness condition is quite technical and typically not trivial
to verify for a specific Hamiltonian system. However, there are some sufficient conditions which imply
steepness, whose verification requires the resolution of algebraic
equalities and inequalities. This motivates the introduction of
the following definition (see \cite{knezevic2008application,Poschel}).

\begin{definition}\label{steepsuff}
Consider the Hamiltonian $h=h(\underline{J})$ for $\underline{J}\in B$ where
$B\subset\real^n$ is an open connected set.
Denote by $\underline{\omega}(\underline{J})$ the gradient of $h$ and by $	\mathcal{Q}(\underline{J})$ its Hessian matrix.
Then:
    \begin{enumerate}
        \item $h(\underline{J})$ is \textit{convex} in $\underline{J}\in B$ if
        \begin{equation*}
        \forall \underline{u}\in\mathbb{R}^n \quad \mathcal{Q}(\underline{J})\underline{u}\cdot\underline{u}=0\Leftrightarrow\underline{u}=\underline{0};
        \end{equation*}
        \item $h(\underline{J})$ is \textit{quasi-convex} in $\underline{J}\in B$ if $\underline{\omega}(\underline{J})\neq0$ and
        \begin{equation*}
        \forall \underline{u}\in\mathbb{R}^n\qquad
        \begin{cases}
        \underline{\omega}(\underline{J})\cdot\underline{u}=0 \\
        \mathcal{Q}(\underline{J})\underline{u}\cdot\underline{u}=0
        \end{cases}
        \Leftrightarrow\underline{u}=\underline{0};
        \end{equation*}
        \item $h(\underline{J})$ is \textit{three-jet non degenerate} in $\underline{J}\in B$ if $\underline{\omega}(\underline{J})\neq0$ and
        \begin{equation*}
        \forall \underline{u}\in\mathbb{R}^n\qquad
        \begin{cases}
        \underline{\omega}(\underline{J})\cdot\underline{u}=0 \\
        \mathcal{Q}(\underline{J})\underline{u}\cdot\underline{u}=0 \\
        \sum_{i,j,k=1}^n\frac{\partial^3 h}{\partial J_i\partial J_j \partial J_k}(\underline{J})u_iu_ju_k=0
        \end{cases}
        \Leftrightarrow\underline{u}=\underline{0}.
        \end{equation*}
    \end{enumerate}
\end{definition}
We remark that the convexity condition is equivalent to require that the Hessian matrix $\mathcal{Q}(\underline{J})$
is positive (or negative) definite in $\underline{J}$. We add also the following definition of
isoenergetically non-degenerate which, for Hamiltonian systems with 2 degrees of freedom, implies
quasi-convexity.

\begin{definition}
    The Hamiltonian $h=h(\underline{J})$ is called \textit{isoenergetically non degenerate} in $\underline{J}\in B$
    with $B\subset \real^n$ open, if
    \begin{equation*}
    \det \begin{pmatrix}
    \frac{\partial^2h}{\partial\underline{J}^2}(\underline{J}) &  \frac{\partial h(\underline{J})}{\partial\underline{J}} \\
    \left(\frac{\partial h(\underline{J})}{\partial\underline{J}}\right)^T & 0
    \end{pmatrix} \neq 0.
    \end{equation*}
\end{definition}

One can prove (see \cite{bambusi2017nekhoroshev}) that, for every
Hamiltonian system with $n$ degrees of freedom, quasi-convexity
implies isoenergetically non-degeneracy: as a consequence, for
two-dimensional Hamiltonian systems, the two conditions are equivalent.

\subsection{Numerical verification of the non-degeneracy conditions}\label{sec: non deg num}
We now apply the above definitions to the Hamiltonian
functions introduced in Section~\ref{model}. We consider the following cases:
\begin{itemize}
    \item the Hamiltonian function related to the $J_2$ problem $\mathcal{H}_{J_2}$, in form of Taylor expansion up to order $15$ in eccentricity and inclination, normalized up to order $12$ with respect to the fast angle $\lambda$;
    we denote the resulting Hamiltonian including the normalized part $\H_{J_2,sec}^{(M)}$ and
    the remainder $\R_{J_2}^{(M)}$ (see Eq. (\ref{hamnew})), as
    \begin{equation*}
    \mathcal{H}_{J_2}(\delta L, P, Q, \lambda, p, q)=\mathcal{H}_{J_2, sec}^{(M)}(\delta L, P, Q, p, q)+\R^{(M)}_{J_2}(\delta L, P, Q, \lambda, p, q)\ .
    \end{equation*}
    Given the practical stability of the semimajor axis established in Section \ref{sec: stabest j2},
    in our computations we set $L=L_*$, i.e., $\delta L=0$;
    \item the Hamiltonian function related to the geolunisolar problem $\mathcal{H}_{gls,sec}$ in
    \equ{hamavfin}, expanded around the forced values of inclination and eccentricity
    (see Section \ref{ssec: forced}) up to order $15$  in eccentricity and inclination, see (\ref{hamavfinaaN}).
    The Hamiltonian $\mathcal{H}_{gls,sec}$ is averaged over the fast angle $\lambda$ and put in
    resonant normal form with respect to the angles $\phi_1$ and $\phi_2$ up to order $12$ in eccentricity and inclination.
    As a consequence, the resulting Hamiltonian $\H_{gls}=\mathcal{H}_{gls, sec}^{(M)}$, including the normalized part
    $\mathcal{Z}_{gls,sec}^{(M)}$, the resonant part $\mathcal{Z}_{gls,res}^{(M)}$ and the remainder $\R_{gls}^{(M)}$ (see Eq. (\ref{Hglsr})),
    has two degrees of freedom and it is the sum of three terms:
    \begin{equation*}
    \mathcal{H}_{gls}(I_1,I_2, \phi_1,\phi_2)=\mathcal{H}_{gls, sec}(I_1,I_2)+\mathcal{H}_{gls,res}(I_1,I_2,\phi_1,\phi_2)+\R_{gls}(I_1,I_2, \phi_1,\phi_2),
    \end{equation*}
    where $\mathcal{H}_{gls, res}$ depends only on the quasi-resonant combination $\phi_1-\phi_2$.
\end{itemize}

To analyze the non-degeneracy conditions, we write the Hamiltonian as the sum of two terms, namely an
integrable Hamiltonian $h$ and a perturbative function $f$.
For the $J_2$-Hamiltonian, we set $h( P, Q)$ to contain all the terms of $\mathcal{H}_{J_2}$ that are independent
on all angles, while the perturbing function $f$ contains all other terms.
For the geolunisolar case, we choose $h(I_1,I_2)$ to be the angle-independent part of the truncation up to order $2$ of
$\mathcal{H}_{gls}$: in this way, the Hessian matrix of $h$ is independent of the actions, and the computations are easier\footnote{We made this particular choice after verifying that, in the chosen domain in the actions, there are no substantial differences between taking all the normalized terms up to order $12$ or only the quadratic truncation.}.

Since the Hamiltonian functions depend on the parameter
$L_*=\sqrt{\mu a_*}$, we select four reference values
for the altitudes that correspond to distances of interest in satellite dynamics:
\begin{itemize}
    \item $3000$ $km$, for near-Earth objects;
    \item $20000$ $km$, for distance of the order of MEO;
    \item $35790$ $km$, for GEO orbits;
    \item $50000$ $km$, for far objects.
\end{itemize}
For each of these values, we check the non-degeneracy conditions of convexity, quasi-convexity and three-jet,
for both the case of the $J_2$-problem and the geolunisolar models in the domain\footnote{From now on, unless otherwise specified, the angles are expressed in $radians$.}
$(e,i)\in \mathcal{D}=[0,0.1]\times[0,0.1]$, which corresponds to a domain in the actions
$\mathcal{D}''=[0,P_{max}]\times[0,Q_{max}]\subset\mathbb{R}^2$, where $P_{max}$, $Q_{max}$ correspond to $e=0.1$, $i=0.1$ and can be computed numerically.

\begin{remark}\label{remark}
We notice that a Hamiltonian $h=h(P,Q)$ (or, equivalently, $h(I_1,I_2)$ in the geolunisolar case) is convex in $\mathcal{D}''\in\mathbb{R}^2$,
if the product of the eigenvalues of the Hessian matrix of $h$ is greater than zero for every $(P,Q)\in \mathcal{D}''$.
Moreover, $h(P,Q)$ is quasi-convex in $\mathcal{D}''\in\mathbb{R}^2$, if for every $(P,Q)\in \mathcal{D}''$ the determinant of the matrix
    \begin{equation}\label{isoe}
    A=  \begin{pmatrix}
    h_{11}(P,Q) & h_{12}(P,Q) & h_1(P,Q) \\
    h_{12}(P,Q) & h_{22}(P,Q) & h_2(P,Q) \\
    h_1(P,Q) & h_2(P,Q) & 0
    \end{pmatrix}
    \end{equation}
    is non zero.
\end{remark}

If the convexity and quasi-convexity tests fail, one can control
the three-jet non-degeneracy condition, that we compute, again, numerically, checking that the system
\begin{equation}\label{threejet}
\begin{cases}
\underline{\omega}(P,Q)\cdot\underline{u}=0 \\
(\partial^2 h(P,Q)\underline{u})\cdot\underline{u}=0 \\
((\partial^3
h(P,Q)\underline{u})\underline{u})\cdot\underline{u}=0
\end{cases}
\end{equation}
evaluated on a grid of values $(P,Q)\in \mathcal{D}''$ admits only
the trivial solution $\underline{u}=(0,0,0)$. Since convexity
implies quasi-convexity and quasi-convexity implies three-jet non-degeneracy,
to identify which of the conditions is satisfied, we
proceed in the following way:
\begin{itemize}
    \item we begin with the convexity test on the product of the eigenvalues: if the product is positive for every value of $(P,Q)\in \mathcal{D}''$, then $h(P,Q)$ is convex;
    \item if the convexity test fails, we pass to the quasi-convexity condition, checking
    the criteria given in Definition~\ref{steepsuff} and Remark~\ref{remark};
    \item if the quasi-convexity test fails, we check the three-jet non-degeneracy through the numerical test based on
    Definition~\ref{steepsuff}.
\end{itemize}

\subsection{Non-degeneracy of the $J_2$ Hamiltonian}\label{ssec: 3j j2}
We start from the convexity test; we denote by $\lambda_1$, $\lambda_2$ the eigenvalues of the
Hessian matrix of $h$.

\begin{table}
    \caption{Values of $\lambda_1\, \lambda_2$ for the $J_2$ model for different altitudes and $(P,Q)\in \mathcal{D}''$.}
    \label{convex j2}
    \begin{center}
        \begin{tabular}{|c|c|}
            \hline
            Altitudes & $\lambda_1\,\lambda_2$ intervals\\
            \hline
            $3000$ $km$ & $[-1.53606\cdot10^{-6},-1.44844\cdot10^{-6}]$ \\
            \hline
            $20000$ $km$ & $[-3.9009\cdot10^{-10}, -3.67893\cdot10^{-10}]$ \\
            \hline
            $35790$ $km$ & $[-3.00586\cdot10^{-16}, -3.19895\cdot10^{-42}]$ \\
            \hline
            $50000$ $km$ & $[-4.3889\cdot10^{-32},-1,30545\cdot10^{-47}]$ \\
            \hline
    \end{tabular}   \end{center}
\end{table}

Table~\ref{convex j2} gives the numerical
values of $\lambda_1\lambda_2$ for different altitudes and $(P,Q)$ in the domain $\mathcal{D}''$
(we recall that, since the values in the Hessian matrix depend on $P$ and $Q$, we have
an interval for $\lambda_1\lambda_2$ instead of a single value).
As one can see, the product of the eigenvalues is always negative \red{or zero within numerical precision level}, leading to the conclusion that the
 Hamiltonian $\H_{J_2}$ is not convex in $\mathcal{D}''$ for the considered altitudes.

\begin{table}
    \caption{Values of $\det{A}$, with $A$ as in (\ref{isoe}) for the $J_2$ model for different altitudes and $(P,Q)\in \mathcal{D}''$.}
    \label{quasiconvex j2}
    \begin{center}
        \begin{tabular}{|c|c|}
            \hline
            Altitudes & $\det{A}$ intervals\\
            \hline
            $3000$ $km$ & $[-1.99418\cdot10^{-10},-1.82748\cdot10^{-10}]$ \\
            \hline
            $20000$ $km$ & $[-2.27575\cdot10^{-15},-2.08271\cdot10^{-15}]$ \\
            \hline
            $35790$ $km$ & $[-1.30574\cdot10^{-17},-1.19784\cdot10^{-17}]$ \\
            \hline
            $50000$ $km$ & $[-5.34517\cdot10^{-19},-4.9035\cdot10^{-19}]$ \\
            \hline
    \end{tabular}   \end{center}
\end{table}

We can then pass to the quasi-convexity test. We consider
the determinant of the matrix $A$ defined in (\ref{isoe}) for
$(P,Q)\in \mathcal{D}''$.
As we can see from Table~\ref{quasiconvex j2}, for every considered altitude the values of $\det{A}$ are equal to zero
within the numerical precision level, leading to the conclusion that the $J_2$ Hamiltonian is not quasi-convex in $\mathcal{D}''$.

The failure of the quasi-convexity for the $J_2$ problem is a relevant fact: as we will see
in Section \ref{ssec: ls nonden}, the effects of the lunisolar attraction will eliminate such degeneracy,
making the total Hamiltonian quasi-convex.

We conclude with the test on the three-jet non-degeneracy condition. To make
the computations quantitative, we solved the system
(\ref{threejet}) for values $(P_i,Q_j)$ on a mesh of $10000$
points in $\mathcal{D}''$. For every pair of values $(P_i,Q_j)$ the
only solution of the system is the trivial one
$\underline{u}=(0,0)$, leading to conclude that the Hamiltonian of the $J_2$ model is three-jet non degenerate in $\mathcal{D}''$.

\subsection{Quasi-convexity of the geolunisolar Hamiltonian}\label{ssec: ls nonden}
As for the $J_2$ model, we start from the convexity test. In this
case, the unperturbed Hamiltonian is a polynomial of degree $2$ in
the actions; then, the Hessian matrix of $h(I_1,I_2)$ does not depend
on the values of $I_1$ and $I_2$, and the same holds for its
eigenvalues. This makes the test on the convexity of the
Hamiltonian easier.

\begin{table}
    \caption{Values of $\lambda_1$ and $\lambda_2$ for different altitudes in the geolunisolar model.}
    \label{convex ls}
    \begin{center}
        \begin{tabular}{|c|c|c|}
            \hline
            Altitudes & $ \lambda_1 $&$ \lambda_2 $\\
            \hline
            $3000$ $km$ & $-11.6416$&$ 3.046$\\
            \hline
            $20000$ $km$ & $-0.185307$&$ 0.0479062$\\
            \hline
            $35790$ $km$ & $-0.0294666$&$ 0.00703346$\\
            \hline
            $50000$ $km$ & $-0.0188881$ &$0.00508227$\\
            \hline
    \end{tabular}   \end{center}
\end{table}
Table \ref{convex ls} shows the values of $\lambda_1$ and $\lambda_2$ for different altitudes. As we can see, in every case the eigenvalues of the Hessian have opposite sign,
showing that the geolunisolar unperturbed Hamiltonian is not convex in $\mathbb{R}^2$, and hence in $\mathcal{D}''$.

As for the quasi-convexity, we check whether the matrix
$A$ defined in (\ref{isoe}) is nondegenerate for every value
$(I_1,I_2)\in \mathcal{D}''$.

\begin{table}
    \caption{Values of $\det{A}$, with $A$ in (\ref{isoe}) in the geolunisolar case for different altitudes and $(I_1,I_2)\in \mathcal{D}''$.}
    \label{quasiconvex ls}
    \begin{center}
        \begin{tabular}{|c|c|}
            \hline
            Altitudes & $\det{A}$ intervals\\
            \hline
            $3000$ $km$ & $[2206.82,2335.21]$ \\
            \hline
            $20000$ $km$ & $[0.0271813,0.0287145]$ \\
            \hline
            $35790$ $km$ & $[0.000309172,0.000323288]$ \\
            \hline
            $50000$ $km$ & $[0.000113523,0.000118622]$ \\
            \hline
    \end{tabular}   \end{center}
\end{table}

From Table \ref{quasiconvex ls} we can see that the determinant of $A$ is
strictly positive for every value of the selected altitudes and every
$(I_1,I_2)\in \mathcal{D}''$. Hence, we conclude that the
Hamiltonian for the geolunisolar case is quasi-convex.
As observed at the end of Section \ref{ssec: 3j j2}, this fact is
highly nontrivial, since it means that the lunisolar perturbation of the $J_2$ model
removes the degeneracy.

\printbibliography

\end{document}